\documentclass[11pt,leqno]{amsart}
\usepackage[T1]{fontenc}
\usepackage[latin1]{inputenc}
\usepackage{amsmath}
\usepackage{amssymb}
\usepackage{latexsym,graphicx,chicago}
\usepackage{graphicx,multirow,rotating}

\newtheorem{thm}{Theorem}[section]
\newtheorem{lem}{Lemma}[section]
\newtheorem{prop}{Proposition}[section]
\newtheorem{cor}{Corollary}[section]

\newtheorem{rem}{Remark}[section]

\def\1{{{\mbox{${\rm{1\negthinspace\negthinspace I}}$}}}}
\newcommand{\eref}[1]{(\ref{#1})}

\newcommand\ind{{ {{1}}\hspace{-0,8mm}{\mathrm I}}}
\newcommand{\esp}{\mathbb{E}}

\newcounter{hypc}

\newcounter{cond}

\begin{document}
\title[Adaptive density estimation for ARCH-type models]{Adaptive density estimation for
general ARCH models}
\author{F. Comte$^{*,1}$}\thanks{$^1$ Universit\'e Paris 5,  MAP5,
UMR CNRS 8145.}
\author{J. Dedecker$^2$}\thanks{$^2$ Universit\'e Paris 6, Laboratoire de
Statistique Th\'eorique et Appliqu\'ee.}
\author{M. L. Taupin $^3$}\thanks{$^3$ IUT de Paris 5 et
Universit\'e Paris Sud, Laboratoire de Probabilit\'es, Statistique
et Mod\'elisation, UMR 8628.}

\begin{abstract}
We consider a model $Y_t=\sigma_t\eta_t$ in which $(\sigma_t)$ is
not independent of the noise process $(\eta_t)$, but $\sigma_t$ is
independent of $\eta_t$ for each $t$. We assume that $(\sigma_t)$ is
stationary and we propose an adaptive estimator of the density of
$\ln(\sigma^2_t)$ based on the observations $Y_t$. Under various dependence structures, the rates of this
nonparametric estimator coincide with the minimax rates obtained
in the i.i.d. case when $(\sigma_t)$ and $(\eta_t)$ are independent,
in all cases where these minimax rates are known. The results apply
to various linear and non linear ARCH processes.
\end{abstract}
\maketitle


\bibliographystyle{chicago}
\section{Introduction}
\setcounter{equation}{0} \setcounter{lem}{0}
\setcounter{thm}{0}

In this paper, we consider the following general ARCH-type model:
$((Y_t, \sigma_t))_{t \geq 0}$ is a strictly stationary sequence of
${\mathbb R}\times {\mathbb R}^+$-valued random variables,
satisfying the equation
\begin{equation}\label{model}
Y_t=\sigma_t\eta_t
\end{equation} where
$(\eta_t)_{t\in {\mathbb Z}}$ is a sequence of  independent and
identically distributed (i.i.d.) random variables with mean zero and
finite variance, and for each $t \geq 0$, the random vector
$(\sigma_i, \eta_{i-1})_{0\leq i \leq t}$ is independent of the
sequence $(\eta_i)_{i \geq t}$.

The model is classically re-written {\it via} a logarithmic
transformation: \begin{equation}\label{convol}
Z_t=X_t+\varepsilon_t, \end{equation} where $Z_t=\ln(Y_t^2)$,
$X_t=\ln(\sigma_t^2)$ and $\varepsilon_t=\ln(\eta_t^2)$. In the
context derived from the model (\ref{model}), $X_t$ and $\varepsilon_t$
are independent for a given $t$, whereas the processes $(X_t)_{t
\geq 0}$ and $(\varepsilon_t)_{t \in {\mathbb Z}}$ are not
independent.

 Our aim is the adaptive estimation of $g$, the
common  distribution of the unobserved variables $X_t= \ln(\sigma_t^2)$,
when the density $f_\varepsilon$ of $\varepsilon_t=\ln(\eta_t^2) $ is
known. More precisely we shall build an estimator  of $g$ without
any prior knowledge on its smoothness, using the observations
$Z_t=\ln(Y_t^2)_t $ and the knowledge of the convolution
kernel $ f_{\varepsilon}$. Since  $X_t$ and $\varepsilon_t$ are
independent for each $t$,  the common density  $f_Z$ of the $Z_t$'s is given by the convolution equation
$f_Z=g *
f_\varepsilon$.

In many papers dealing with ARCH models, $\varepsilon_t$ is assumed to be Gaussian or the log of a squared
Gaussian (when $\eta_t$ is Gaussian, see van Es {\it et
al.}~\citeyear{vanes05} or in slightly different contexts van Es
{\it et al.}~\citeyear{vanes03}, Comte and
Genon-Catalot~\citeyear{CGC}). Our setting is more general since we
consider various type of error densities. More precisely, we assume that $f_\varepsilon$
belongs to some class of smooth functions described below:
 there exist nonnegative numbers
$\kappa_0$, $\gamma$, $\mu$, and $\delta$  such
that the fourier transform $f_{\varepsilon}^*$ of $f_\varepsilon$ satisfies
\begin{equation}
\kappa_0(x^2+1)^{-\gamma/ 2}\exp\{-\mu\vert x\vert^\delta\}\leq
|f_\varepsilon^*(x)| \leq 
\kappa_0'(x^2+1)^{-\gamma/ 2}\exp\{-\mu\vert x\vert^\delta\}.
\label{condfeps}
\end{equation}
Since $f_\varepsilon$ is known, the constants $\mu, \delta,
\kappa_0,$ and $\gamma$  defined in \eref{condfeps} are known. When
$\delta=0$ in \eref{condfeps}, the errors are called ``ordinary
smooth'' errors. When $\mu>0$ and $\delta>0$, they are called
``super smooth''. The standard examples for super smooth densities
are Gaussian or Cauchy distributions (super smooth of order
$\gamma=0, \delta=2$ and $\gamma=0, \delta=1$ respectively). When
$\varepsilon_t=\ln(\eta_t^2)$ with $\eta_t\sim {\mathcal N}(0,1)$ as
in van Es {\it et al.}~\citeyear{vanes03,vanes05}, then $\varepsilon_t$ is super-smooth
with $\delta=1, \gamma=0$ and $\mu= \pi/2$. An example of ordinary
smooth density is the Laplace distribution, for which $\delta=\mu=0$
and $\gamma=2$.

In density deconvolution of i.i.d variables the  $X_t$'s and the
$\varepsilon_t$'s are i.i.d. and the  sequences $(X_t)_{t \geq 0}$ and
$(\varepsilon_t)_{t \in {\mathbb Z}}$ are independent (for short we shall refer to this case as the i.i.d. case).
In the setting of Model (\ref{convol}), the classical assumptions of
independence between the processes $(X_t)_{t \geq 0}$ and
$(\varepsilon_t)_{t \in {\mathbb Z}}$ are no longer satisfied and the
tools for deconvolution have to  be revisited.

As in density deconvolution for i.i.d. variables, the slowest rates of
convergence for estimating $g$ are obtained for super smooth error
densities. For instance, in the i.i.d case, when $\varepsilon_t$ is Gaussian or the log
of a squared Gaussian and $g$ belongs to some  Sobolev class, the
minimax rates are negative powers of $\ln(n)$ (see
Fan~\citeyear{Fan91}). Nevertheless, it has been noticed by several
authors (see Pensky and Vidakovic \citeyear{Penskyvidakovic99},
Butucea~\citeyear{butucea}, Butucea and Tsybakov \citeyear{ButTsyb},
Comte {\it et al.}~\citeyear{CRT1}) that the rates are improved if
$g$ has stronger smoothness properties.
So, we describe the smoothness
properties of $g$ by the set
\begin{equation}
\label{super} \mathcal{S}_{s,r,b}(C_1)=\Big \{\psi\; \mbox{ such that }
\int_{-\infty}^{+\infty} |\psi^*(x)|^2(x^2+1)^{s}\exp\{2b |x|^{r} \}
dx\leq C_1 \Big \}
\end{equation}
for $s,r,b$ unknown non negative numbers. When $r=0$,
the class $\mathcal{S}_{s,r,b}(C_1)$ corresponds to a Sobolev ball. When $r>0, b>0 $
functions belonging to $\mathcal{S}_{s,r,b}(C_1)$ are infinitely many times differentiable.

Our estimator of $g$ is constructed by   minimizing an appropriate
penalized contrast function only depending on the observations and
on $f_\varepsilon$.  It is chosen in a purely data-driven way among a collection of
non-adaptive estimators. We start by the study of those non-adaptive estimators and show
that their mean integrated squared  error (MISE) has the same
order as in the i.i.d. case. In particular they reach the
minimax rates of the i.i.d. case in all cases where they are known (see
Fan~\citeyear{Fan91}, Butucea~\citeyear{butucea} and Butucea and
Tsybakov \citeyear{ButTsyb}). Next we prove that the MISE of our adaptive estimator is of the same
order as the MISE of the best non-adaptive estimator, up to  some possible
negligible logarithmic loss in one case.

In their 2005 paper, van Es \textit{et al.}~\citeyear{vanes05} have considered the case  where $\eta_t$ is Gaussian, the density $g$ of $X_t$ is twice  differentiable, and the process $(Z_t,X_t)$ is $\alpha$-mixing. Here  we consider various types of error density, and we do not make any  assumption on the smoothness of $g$: this is the advantage of the adaptive procedure.
We shall consider two types of dependence properties, which are satisfied by many ARCH processes.
First we shall use the classical $\beta$-mixing properties of
general ARCH models, as recalled in Doukhan~\citeyear{Doukhan} and described in more details
in Carrasco and Chen~\citeyear{CCET}.
But we also illustrate that new recent coefficients can be
used in our context, which allow an easy characterization of the
dependence properties in function of the parameters of the models. Those new
dependence coefficients, recently defined and studied in Dedecker and
Prieur~\citeyear{JDCPCoeff}, are interesting and powerful because
they require much lighter conditions on the models. Such ideas
have been popularized by  Ango Nz\'e and
Doukhan~\citeyear{ND} and Doukhan {\it et al.} \citeyear{DTW}. For instance, these coefficients  allow to deal with the
general ARCH$(\infty)$ processes defined by Giraitis {\it et
al.}~\citeyear{GKL}.

The paper  is organized as follows. Many  examples are described in
Section \ref{depend}, together with their dependence properties. The
estimator is defined in Section \ref{estim}.  The MISE
bounds are given in Section \ref{risk}, and the proofs are given in Section 5.

\section{The model and its dependence properties}\label{depend}
\setcounter{equation}{0}
\setcounter{lem}{0}
\setcounter{thm}{0}

\subsection{Models and examples}
A particular case of model (\ref{model}) is
\begin{equation}\label{shift}
Y_t=\sigma_t \eta_t, \ \text{with} \   \sigma_t= f( \eta_{t-1},
\eta_{t-2}, \ldots)
\end{equation}
for some measurable function $f$. Another important case is
\begin{equation}\label{markov}
Y_t= \sigma_t \eta_t, \ \text{with} \ \sigma_t= f(\sigma_{t-1},
\eta_{t-1}) \ \text{and $\sigma_0$ independent of $(\eta_t)_{t\geq
0}$,}
\end{equation}
 that is $\sigma_t$ is a stationary  Markov chain.

We begin with models satisfying a recursive
equation,  whose stationary solution satisfies (\ref{shift}).
The original ARCH model as introduced by Engle~\citeyear{engle} was given by
\begin{equation}\label{arch1}
 Y_t= \sqrt{a+b
Y_{t-1}^2}\eta_{t}, \  a\geq 0, b \geq  0
\end{equation} It has been
generalized by Bollerslev~\citeyear{bollerslev} with the class of
GARCH$(p,q)$ models defined by $Y_t=\sigma_t \eta_t$ and
\begin{equation}\label{garchrel} \sigma_t^2= a +\sum_{i=1}^p
a_iY_{t-i}^2 +\sum_{j=1}^q b_j \sigma_{t-j}^2\end{equation} where
the coefficients $a, a_i, i=1, \dots, p$ and $b_j, j=1, \dots, q$
are all positive real numbers. Those processes were studied from the
point of view of existence and stationarity of solutions by Bougerol
and Picard~\citeyear{BP2,BP1} and Ango Nz\'e~\citeyear{angonze}.
Under the condition  $\sum_{i=1}^p a_i + \sum_{j=1}^q b_j <1$, this
model has a unique stationary solution of the form (\ref{shift}).

Many extensions have been proposed since then. A general linear
example of model  is given by the ARCH$(\infty)$ model described by
Giraitis {\it et al.}~\citeyear{GKL}:
\begin{equation}\label{archinfty}
 \sigma^2_t = a+ \sum_{j=1}^{\infty}
a_j Y^2_{t-j},\end{equation} where $a\geq 0$ and $a_j\geq 0$. Again
if $ \sum_{j \geq 1} a_j < 1$, then there exists a unique  strictly
stationary solution to (\ref{archinfty}) of the form (\ref{shift}).

For the models satisfying (\ref{markov}), let us cite first  the
so-called augmented GARCH$(1,1)$ models introduced by
Duan~\citeyear{DUAN}:
\begin{equation}\label{duan}
\Lambda(\sigma^2_t)=c(\eta_{t-1})\Lambda(\sigma^2_{t-1}) +
h(\eta_{t-1}),\end{equation} where $\Lambda$ is an increasing and
continuous function on ${\mathbb R}^+$.  We refer to
Duan~\citeyear{DUAN} for numerous examples of more standard models
belonging to this class. There exists a stationary solution  to
(\ref{duan}), provided $c$ satisfies the condition $A_2^*$ given in
Carrasco and Chen~\citeyear{CCET} (this condition is satisfied as
soon as ${\mathbb E}(|c(\eta_0)|^s)< 1$ and ${\mathbb
E}(|h(\eta_0)|^s)< \infty$ for integer $s\geq 1$, see the
condition $A_2$ of the same paper).
An example of the model (\ref{duan}) is the threshold ARCH model
(see Zako\"{\i}an~\citeyear{zakoian}):
\begin{equation}\label{tarch} \sigma_{t}= a+b\sigma_{t-1}\eta_{t-1}
\1_{\{\eta_{t-1}>0\}} -c \sigma_{t-1}
\eta_{t-1}\1_{\{\eta_{t-1}<0\}}, \ a,b,c>0\end{equation} for which
$c(\eta_{t-1})=b\eta_{t-1} \1_{\{\eta_{t-1}>0\}} -c
\eta_{t-1}\1_{\{\eta_{t-1}<0\}}$ and $h=a$. In particular, the
condition for the stationarity is satisfied as soon as $b\vee c <1$.

Other  models satisfying (\ref{markov}) are the  non linear ARCH
models (see Doukhan~\citeyear{Doukhan}, p. 106-107), for which:
\begin{equation}\label{model1}
 \sigma_t=f(\sigma_{t-1}\eta_{t-1}).
\end{equation}
There exists a stationary solution to (\ref{model1}) provided that
the density of $\eta_0$ is  positive on a neighborhood of 0 and
$\lim\sup_{|x|\rightarrow \infty} |f(x)/x| <1$.

In the next section, we define the dependence coefficients that we
shall use in this paper, and we give the dependence properties of
the models (\ref{arch1})-(\ref{model1}) in terms of these
coefficients.

\subsection{Measures of dependence}
\label{couplage}

Let $(\Omega, {\mathcal A}, {\mathbb P})$ be a probability space.
Let $W$ be a random vector with values in a Banach space
$({\mathbb B}, \|\cdot \|_{\mathbb B})$, and let ${\mathcal M}$ be
a $\sigma$-algebra of ${\mathcal A}$. Let ${\mathbb
P}_{W|{\mathcal M}}$ be a conditional distribution of $W$ given
${\mathcal M}$, and let $P_W$ be the distribution of $W$. Let
${\mathcal B}({\mathbb B})$ be the Borel $\sigma$-algebra on
$({\mathbb B}, \| \cdot \|_{\mathbb B})$, and let
$\Lambda_1({\mathbb B})$ be the set of 1-Lipschitz functions from
$({\mathbb B}, \|\cdot\|_{\mathbb B})$ to ${\mathbb R}$. Define
now
\begin{eqnarray*}
\beta({\mathcal M}, \sigma(W))&=& {\mathbb E}\Big (\sup_{A \in
{\mathcal B}({\mathcal X})} |{\mathbb P}_{W|{\mathcal
M}}(A)-{\mathbb P}_W(A)|
\Big ) \, , \\
\text{and if ${\mathbb E}(\|W\|_{\mathbb B}) < \infty$,} \quad \tau
({\mathcal M}, W)&=& {\mathbb E}\Big (\sup_{f \in
{\Lambda_1}({\mathbb B})} |{\mathbb P}_{W|{\mathcal M}}(f)-{\mathbb
P}_W(f)| \Big ) \, .
\end{eqnarray*}
The coefficient $\beta({\mathcal M}, \sigma(W))$ is the usual
mixing coefficient, introduced by Rozanov and Volkonskii
\citeyear{VolkRozan}. The coefficient $\tau({\mathcal M}, W)$ has
been introduced by Dedecker and Prieur \citeyear{JDCPCoeff}.

Let $(W_t)_{t \geq 0}$ be a strictly stationary sequence of
${\mathbb R}^2$-valued random variables. On ${\mathbb R}^2$, we put
the norm $\|x-y\|_{\mathbb R^2}=|x_1-y_1|+ |x_2-y_2|$. For any $k
\geq 0$, define the coefficients
\begin{equation}\label{defbeta1}
\beta_{1}(k)=  \beta(\sigma (W_0), \sigma(W_{k})), \quad \text{ and
if } {\mathbb E}(\|W_0\|_{{\mathbb R}^2}) < \infty, \quad
\tau_{1}(k)= \tau(\sigma(W_0), W_{k}).
\end{equation}

On $({\mathbb R}^2)^l$, we put the norm $\|x-y\|_{({\mathbb
R^2})^l}=l^{-1}(\|x_1-y_1\|_{\mathbb R^2}+ \cdots
+\|x_l-y_l\|_{\mathbb R^2})$. Let ${\mathcal M}_i=\sigma(W_k, 0\leq
k \leq i)$. The coefficients $\beta_{\infty}(k)$ and
$\tau_{\infty}(k)$ are defined by
\begin{equation}\label{betainfini}
\beta_{\infty}(k)= \sup_{i \geq 0}\sup_{ l \geq 1}  \left \{
\beta({\mathcal M}_i, \sigma(W_{i_1}, \ldots , W_{i_l})), i+ k \leq
i_1 < \cdots < i_l \right \},
\end{equation}
and if ${\mathbb E}(\|W_1\|_{\mathbb R^2}) < \infty$,
\begin{equation} \label{tauinfini}
\tau_{\infty}(k)= \sup_{i \geq 0} \sup_{ l \geq 1}  \left \{
 \tau({\mathcal M}_i, (W_{i_1}, \ldots , W_{i_l})), i+k \leq i_1 < \cdots < i_l \right \}.
\end{equation}

 We say that the process $(W_t)_{t \geq 0}$ is $\beta$-mixing (resp.
$\tau$-dependent) if the coefficients $\beta_{\infty}(k)$ (resp.
$\tau_{\infty}(k)$) tend to zero as $k$ tends to infinity. We say
that it is geometrically  $\beta$-mixing (resp. $\tau$-dependent), if
there exist $ a >1$ and $C>0$ such that
 $\beta_{\infty}(k)\leq C a^k$ (resp. $\tau_{\infty}(k)\leq C a^k$) for all $k \geq 1$.

We now recall the coupling properties associated with the dependency
coefficients. Assume that $\Omega$ is rich enough, which means
that there exists $U$ uniformly distributed over $[0,1]$ and
independent of ${\mathcal M} \vee \sigma (W)$. There exist two
${\mathcal M} \vee \sigma(U) \vee \sigma (W)$-measurable random
variables $W_1^\star$ and $W_2^\star$ distributed as $W$ and
independent of ${\mathcal M}$ such that
\begin{equation}\label{cber}
\beta ({\mathcal M}, \sigma(W))= {\mathbb P}(W \neq W_1^\star) \quad
\text{and} \quad  \tau({\mathcal M}, W)= {\mathbb E}(\|W - W_2^\star\|_{\mathbb B}) \, .
\end{equation}
The first equality in (\ref{cber}) is due to Berbee
\citeyear{Berbee}, and the second one has been established in
Dedecker and Prieur \citeyear{JDCPCoeff}, Section 7.1.\\

 As consequences of the coupling properties \eref{cber}, we have the following covariance inequalities.
 Let $\| \cdot
\|_{\infty, {\mathbb P}}$ be the ${\mathbb L}^\infty (\Omega,
{\mathbb P})$-norm.
For two measurable functions $f, h$ from ${\mathbb R}$
to ${\mathbb C}$, we have
\begin{equation}\label{ibrb}
|\text{Cov}(f(Y), h(X))| \leq 2 \|f(Y)\|_{\infty, {\mathbb P}}
\|h(X)\|_{\infty, {\mathbb P}} \,  \beta (\sigma (X), \sigma (Y))
\, .
\end{equation}
 Moreover, if $\text {Lip}(h)$ is the Lipschitz coefficient of $h$,
\begin{equation}\label{evident}
 \quad |\text{Cov }(f(Y), h(X))|
\leq \|f(Y)\|_{\infty, {\mathbb P}} \text{Lip}(h) \, \tau (\sigma
(Y), X) \, .
\end{equation}
Thus, using that $t \rightarrow e^{ixt}$ is
$|x|$-Lipschitz, we obtain the  bounds
\begin{eqnarray}
\label{covineg}
      |\mbox{Cov}(e^{ixZ_1},e^{ixX_k})| \leq 2 \beta_{1} (k-1)
    \quad \text{and }
\quad     |\mbox{Cov}(e^{ixZ_1},e^{ixX_k})| &\leq&   |x| \tau_{1}
   (k-1).
\end{eqnarray}

 \subsection{Application to ARCH models}

For the models (\ref{model}) and (\ref{convol}), the $\beta$-mixing
coefficients of the process \begin{eqnarray}
(W_t)_{t\in \mathbb{Z}}=((Z_t, X_t))_{t \in {\mathbb Z}}\label{Wt}
\end{eqnarray} are
smaller than that of   $((Y_t,\sigma_t))_{t \in {\mathbb Z}}$
(because of the inclusion of $\sigma$-algebras). If we assume that
in all cases the $\eta_t$'s are centered with unit variance and
admit a density with respect to the Lebesgue measure,
then
\begin{itemize}
\item The process $((Y_t,\sigma_t))_{t \in {\mathbb Z}}$ defined by Model (\ref{arch1})
is  geometrically $\beta$-mixing  as soon as  $0<b<1$.
\item The process $((Y_t,\sigma_t))_{t \in {\mathbb Z}}$ defined by Model
(\ref{garchrel}) is  geometrically $\beta$-mixing, as soon as
$\sum_{i=1}^p a_i+\sum_{j=1}^q b_j<1$ (see Carrasco and
Chen~\citeyear{CCCRAS,CCET}).
\item The process $((Y_t,\sigma_t))_{ t \in {\mathbb Z}}$ defined
by Model (\ref{duan}) is geometrically $\beta$-mixing as soon as:
the density of $\eta_0$ is   positive on an open set containing 0;
$c$ and $h$ are polynomial functions; there exists an integer $s\geq
1$ such that $|c(0)|<1$, ${\mathbb E}(|c(\eta_0)|^s)<1$, and
${\mathbb E}(|h(\eta_0)|^s) <\infty$. See Proposition 5 in Carrasco
and Chen~\citeyear{CCET}.
\item The process $((Y_t,\sigma_t))_{t \in {\mathbb Z}}$ defined by Model (\ref{tarch}) is
 geometrically $\beta$-mixing  as soon as  $0<b\vee c<1$.
\item The process $((Y_t,\sigma_t))_{t \in {\mathbb Z}}$ defined by Model (\ref{model1}) is
 geometrically $\beta$-mixing as soon as the density of $\eta_0$ is positive on a neighborhood of $0$ and $\lim\sup_{|x|\rightarrow
+\infty} |f(x)/x| <1$ (see Doukhan~\citeyear{Doukhan}, Proposition 6
page 107).
\end{itemize}
 Note that some other extensions to nonlinear models having
stationarity and dependency properties can be found in Lee and
Shin~\citeyear{LS}.

Concerning the $\tau$-dependence, here is a general method to handle
the models (\ref{shift}) and (\ref{markov}). The following
Proposition will be proved in appendix (see Ango Nz\'e and Doukhan~\citeyear{ND} and Doukhan {\it et al.}~\citeyear{DTW} for related results).

\begin{prop}\label{taumix} Let $Y_t$ and $\sigma_t$ satisfy either (\ref{shift}) or (\ref{markov}).
For Model (\ref{shift}), let $(\eta'_t)_{t \in {\mathbb Z}}$ be
an independent copy of $(\eta_t)_{t \in {\mathbb Z}}$, and for
$t>0$, let $\sigma_t^*= f(\eta_{t-1}, \ldots, \eta_1, \eta'_0,
\eta'_{-1}, \ldots)$. For Model (\ref{markov}), let $\sigma_0^*$
be a
 copy of $\sigma_0$ independent of $(\sigma_0, \eta_t)_{t \in {\mathbb Z}}$,
 and for $t>0$
let $\sigma_t^*= f(\sigma_{t-1}^*, \eta_{t-1})$. 
  Let $\delta_n$ be a non increasing
sequence such that
\begin{equation}\label{(2.10)}
2{\mathbb E}(|\sigma_n^2- (\sigma_n^*)^2|) \leq \delta_n \, .
\end{equation}
Then
\begin{enumerate}
\item The process $((Y_t^2, \sigma_t^2))_{t \geq 0}$ is $\tau$-dependent with $\tau_{\infty}(n) \leq \delta_n$.
\item
Assume  that $Y_0^2$, $\sigma_0^2$ have densities satisfying
$\max(f_{\sigma^2}(x), f_{Y^2}(x))\leq C |\ln(x)|^\alpha x^{-\rho}$
in a neighborhood of $0$, for some $\alpha\geq 0$ and $ 0\leq \rho
<1$. The process $((X_t, Z_t))_{t \geq 0}$ is $\tau$-dependent with
$\tau_{\infty}(n) =O((\delta_n)^{(1- \rho)/(2 - \rho)}
|\ln(\delta_n)|^{(1+\alpha)/(2-\rho)})$.
\end{enumerate} Consider Model (\ref{archinfty}), and assume that $c= \sum_{j \geq 1} a_j < 1 $. Let then $((Y_t,
\sigma_t))_{t \in {\mathbb Z}}$ be the unique strictly stationary
solution of the form (\ref{shift}). Then (\ref{(2.10)}) holds with
$$
\delta_n= O \Big ( \inf_{1 \leq k \leq n} \Big \{ c^{n/k} +
\sum_{i=k+1}^\infty a_i \Big \} \Big )\, .
$$
\end{prop}
Note that if $\sigma_0^2$ and $\eta_0^2$ have bounded densities,
then $f_{Y^2}(x))\leq C |\ln(x)|$ in a neighborhood of 0, so that
Proposition \ref{taumix}(2) holds with $\rho=0$ and $\alpha=1$.

Under the assumptions of Proposition \ref{taumix}(2), we obtain for
Model (\ref{archinfty}) the following rates for $((X_t, Z_t))_{t
\geq 0}$:
\begin{itemize}
\item If $a_j=0$, for $j \geq J$, then $((X_t,Z_t))_{t \geq 0}$ is geometrically
$\tau$-dependent.
\item If $a_j=O(b^j)$ for some $b<1$ then $\tau_{\infty}(n)
=O(\kappa^{\sqrt n})$ for some $\kappa<1$.
\item If $a_j=O(j^{-b})$ for some $b>1$ then
$\tau_{\infty}(n) =O(n^{-b(1- \rho)/(2- \rho)}
(\ln(n))^{(b+2)(1+\alpha)/2})$.
\end{itemize}
 For more general models than (\ref{archinfty}), we refer to
Doukhan {\it et al.} \citeyear{DTW}.

For Model (\ref{markov}), if there exists $\kappa<1$ such that
\begin{equation}\label{lip}
   {\mathbb E}(|(f(x, \eta_0))^2-(f(y, \eta_0))^2|) \leq \kappa |x^2-y^2|\, ,
\end{equation}
then one can take $\delta_n= 4{\mathbb E}(\sigma_0^2) \kappa^n $.
Hence, under the assumptions  of Proposition \ref{taumix}(2),
$((X_t, Z_t))_{t
>0}$ is geometrically $\tau$ dependent. An example of Markov chain
satisfying (\ref{lip}) is the autoregressive model $\sigma_t^2=
h(\sigma_{t-1}^2)+ r(\eta_{t-1})$ for some $\kappa$-lipschitz
function $h$.

\section{The estimators}\label{estim}
\setcounter{equation}{0} \setcounter{lem}{0} \setcounter{thm}{0}
For two complex-valued functions
$u$ and $v$ in $\mathbb{L}_2(\mathbb{R})\cap
\mathbb{L}_1(\mathbb{R})$, let $ u^*(x)=\int e^{itx}u(t)dt$,
$u*v(x)=\int u(y)v(x-y)dy,$ and $\langle u,v\rangle =\int
u(x)\overline{v}(x)dx$ with $\overline{z}$ the conjugate of a
complex number $z$. We also denote by $\|u\|_1=\int |u(x)|dx$,
$\|u\|^2=\int |u(x)|^2 dx$, and $\|u\|_\infty=\sup_{x \in
\mathbb{R}}|u(x)|.$

\subsection{The projection spaces }

Let $\varphi(x)=\sin(\pi x)/(\pi x)$. For $m\in {\mathbb N}$ and
$j\in {\mathbb Z}$, set $\varphi_{m,j}(x) =\sqrt{m}
\varphi(mx-j).$ The functions $\{\varphi_{m,j}\}_{j \in
\mathbb{Z}}$ constitute an orthonormal system in ${\mathbb
L}^2({\mathbb R})$ (see e.g.  Meyer \citeyear{MeyerI}, p.22). Let
us define
$$S_m= \overline{{\rm span}}\{\varphi_{_{m,j}}, \;
j\in \mathbb{Z}\}, m\in {\mathbb N}.$$ The space $S_m$ is exactly
the  subspace of ${\mathbb L}_2({\mathbb R})$ of functions having
a Fourier transform with compact support contained in $[-\pi m,
\pi m]$.
The orthogonal projection of $g$ on $S_m$ is $g_m=\sum_{j\in
{\mathbb Z}} a_{m,j}(g) \varphi_{m,j}$ where $a_{m,j}(g) =
<\varphi_{m,j},g>$. To obtain representations having a finite
number of  "coordinates", we introduce
$$S_{m}^{(n)}= \overline{\rm span }\left\{\varphi_{m,j}, |j|\leq
k_n\right\}$$ with integers
$k_n$ to be specified later. The family
$\{\varphi_{m,j}\}_{\vert j\vert \leq k_n}$ is an orthonormal
basis of $S_m^{(n)}$ and the orthogonal projections of $g$ on $S_{m}^{(n)}$ is
 given by $g_{m}^{(n)}=\sum_{|j|\leq k_n}
a_{m,j}(g) \varphi_{m,j}$. Subsequently a space $S_m^{(n)}$ will
be referred to as a "model" as well as a "projection space".

\subsection{Construction of the minimum contrast estimators}
\label{nonpen}

We subsequently assume  that \begin{equation}\label{feps}
f_{\varepsilon} \mbox{  belongs to } {\mathbb L}_2({\mathbb R})
\mbox{ and is such that }  \forall x\in
\mathbb{R},\,f_\varepsilon^*(x)\not=0.\end{equation}  Note that
the square integrability of $f_{\varepsilon}$ and \eref{condfeps}
require that $\gamma> 1/2$ when $\delta=0$. Under Condition
(\ref{feps})  and for or $t$ in $S_m^{(n)}$,  we define the contrast function

\begin{equation*}\label{u*}
\gamma_n(t)=\frac{1}{n}\sum_{i=1}^n \left[\|t\|^2
-2u_t^*(Z_i)\right], \;\;\; \mbox{ with } \;\;\;u_t(x) = \frac 1{2
\pi} \left(\frac{t^*(-x)}{f_{\varepsilon}^*(x)}\right).
\end{equation*}
Then, for an arbitrary fixed integer $m$, an estimator of $g$
belonging to $S_m^{(n)}$ is defined by
\begin{equation}\label{tronque}  \hat g_m^{(n)} = \arg\min_{t\in S_m^{(n)}}
\gamma_n(t).\end{equation}
By using Parseval and inverse Fourier formulae we obtain  that
$\mathbb{E}\left[u_t^*(Z_i)\right] =\langle t, g\rangle,$ so that
$\mathbb{E}(\gamma_n(t))=\|t-g\|^2 -\|g\|^2 $ is minimal when $t=g$.
This shows that $\gamma_n(t)$ suits well for the estimation of $g$.
It is easy to see that $$\hat g_m^{(n)} = \sum_{\vert j\vert\leq
k_n} \hat a_{m,j} \varphi_{m,j} \;\; \mbox{with}\;\; \hat a_{m,j}=
\frac{1}{n} \sum_{i=1}^n u_{\varphi_{m,j}}^*(Z_i), \text{ and }
\mathbb{E}(\hat a_{m,j})= <g,\varphi_{m,j}>=a_{m,j}(g).$$

\subsection{Minimum penalized contrast estimator}

The minimum penalized estimator
of $g$ is defined  as $\tilde g=\hat g_{\hat m_g}^{(n)}$ where $\hat
m_g$ is chosen in a purely data-driven way. The main point of the
estimation procedure lies in the choice of $m=\hat m$ (or equivalently in the choice of model
$S_{\hat m}^{(n)}$) involved in the
estimators $\hat g_m^{(n)}$ given by (\ref{tronque}), in order to mimic
the oracle parameter
\begin{eqnarray}
\label{oraclepar}
\breve m_g=\arg\min_{m}\mathbb{E}\parallel \hat  g_m^{(n)}-g\parallel_2^2.
\end{eqnarray}
The model selection is performed in an automatic way, using the
following penalized criteria
\begin{equation}\label{estitronc}
\tilde g=\hat g^{(n)}_{\hat m} \mbox{ with } \hat m= \arg\min_{m\in
\{1, \cdots , m_n \}} \left[\gamma_n(\hat g_m^{(n)}) + \; {\rm
pen}(m)\right],
\end{equation}
where $\mbox{pen}(m)$ is a penalty function that
depends on $f_\varepsilon^*(\cdot)$ through $\Delta(m)$ defined by
\begin{eqnarray}
\Delta(m)=\frac{1}{2\pi}\int_{-\pi m}^{\pi m}\frac{1}{\vert
f_\varepsilon^*( x)\vert^2}dx.\label{Delta1}
\end{eqnarray}
The key point in the dependent context is to find a penalty
function not depending on the dependency coefficients such that
$$\mathbb{E}\parallel \tilde g -g\parallel^2\leq C\inf_{m\in
\{1,\cdots,m_n\}}\mathbb{E}\parallel \hat g_m^{(n)}-g\parallel^2.$$
In that way,  the  estimator $\tilde g$ is adaptive since it achieves the best rate among the estimators
$\hat g_m^{(n)}$, without any prior knowledge on the smoothness on $g$.

\section{Density estimation bounds}\label{risk}
\setcounter{equation}{0}
\setcounter{lem}{0}
\setcounter{thm}{0}

>From now on, the dependence coefficients are defined as in \eref{defbeta1},
\eref{betainfini} and \eref{tauinfini}
with $(W_t)_{t \in \mathbb{Z}}=((Z_t, X_t))_{t \in\mathbb{Z}}$.

\subsection{Rates of convergence of the minimum contrast estimators $\hat g_m^{(n)}$}

Subsequently, the density $g$ is assumed to satisfy the following
assumption:
\begin{equation}\label{momg} g\in \mathbb{L}_2(\mathbb{R}), \mbox{and there exists } M_2>0, \; \int x^2g^2(x)dx \leq M_2<\infty.
\end{equation}
Assumption \eref{momg}, which is due to the construction of the
estimator, already appears in density deconvolution in the  independent framework in
Comte {\it et al.}~\citeyear{CRT2,CRT1}. It is important to note that Assumption \eref{momg} is
very unrestrictive.  In particular, all densities having tails of order $|x|^{-(s+1)}$ as $x$ tends to
infinity satisfy \eref{momg}  only if $s>1/2$. One can cite for
instance the Cauchy distribution or all stable distributions with
exponent $r>1/2$ (see Devroye~\citeyear{Devroye86}). The L\'evy
distribution, with exponent $r=1/2$ does not satisfies \eref{momg}.

Note that (\ref{momg}) is fulfilled if $g$ is bounded by $M_0$ and
${\mathbb E}(X_1^2)\leq M_1<+\infty$, with $M_2=M_0M_1$.

The order of the MISE of $\hat g_m^{(n)}$ is given in the following proposition.

\begin{prop}
\label{Vitss} If  \eref{feps} and \eref{momg} hold, then $\hat g_m^{(n)}$ defined
by \eref{tronque} satisfies
$$ {\mathbb E}\|g-\hat g_m^{(n)}\|^2\leq
\|g-g_m\|^2+ \frac{m^2(M_2+1)}{k_n} + \frac{2\Delta(m)}{n}+
\frac{2R_{m}}{n},
$$
where \begin{equation}\label{Rm}
 R_{m}=\frac{1}{\pi}\sum_{k=2}^{n}\int_{-\pi
m}^{\pi m}\Big\vert \frac{ {\mathrm
{Cov}}\left(e^{ixZ_1},e^{ixX_k}\right)}{f_{\varepsilon}^*(-x)}\Big\vert
dx.
\end{equation}
Moreover, $R_m \leq \min ( R_{m, \beta},  R_{m,
\tau})$, where
\begin{equation*}
  R_{m, \beta}= 4 \Delta_{1/2}(m) \sum_{k=1}^{n-1}\beta_{1}(k)
\quad \text{and} \quad
  R_{m, \tau} = \pi m \Delta_{1/2}(m) \sum_{k=1}^{n-1}\tau_{1}(k) \, ,
\end{equation*}
with $\beta_1$, $\tau_1$ defined by (\ref{defbeta1}), and where
\begin{equation}\label{midelta}
\Delta_{1/2}(m)=\frac{1}{2\pi}\int_{-\pi m}^{\pi m}\frac{1}{\vert
f_\varepsilon^*( x)\vert}dx.
\end{equation}
\end{prop}
This proposition requires several comments.

As usual, the order of the risk is given by a bias term $\parallel
g_m-g\parallel^2+m^2(M_2+1)/k_n$ and a variance term $2\Delta(m)/n+2R_m/n$.
 As in density
deconvolution for i.i.d. variables, the variance term $2\Delta(m)/n+2R_m/n$ depends on the rate of decay
of the Fourier transform of $f_\varepsilon$. It is the sum of the variance term appearing in density
deconvolution for i.i.d. variables
$2\Delta(m)/n$ and of an additional term $2R_m/n$. This last term $R_m$ involves
the dependency
coefficients and the quantity $\Delta_{1/2}(m)$, which is specific to the ARCH problem.
The point is that, as in the i.i.d. case, the main order term in the variance part is $\Delta(m)/n$,
which does not involve the dependency coefficients. In other words, the
dependency coefficients only appear in front of the additional and negligible
term $\Delta_{1/2}(m)/n$, specific to ARCH models.

The bias term is the sum of the usual bias term $\parallel
g_m-g\parallel^2$, depending on
the smoothness properties of $g$,  and on an additional term $m^2(M_2+1)/k_n$. With a suitable
choice of $k_n$, not depending on $g$, this last term  is negligible with respect to
the variance term.

Concerning the main variance term, $\Delta(m)$ given by (\ref{Delta1}) has the same order as 
$$\Gamma(m)= (1+(\pi m)^2)^{\gamma}(\pi m)^{1-\delta}
\exp\left\{2\mu(\pi m)^{\delta}\right\},$$ up to some constant bounded by \begin{eqnarray}
\label{lambda1}
\qquad\quad\lambda_1(f_\varepsilon,\kappa_0)=\frac{1}
{\kappa_0^{2}\pi R(\mu,\delta)}, \mbox{ where }\, R(\mu,\delta)=
\mbox{1}\!\mbox{I}_{\{\delta=0\}} + 2\mu\delta
\mbox{1}\!\mbox{I}_{\{\delta>0\}} .
\end{eqnarray}

 The rates resulting from Proposition
\ref{Vitss} under \eref{condfeps} and \eref{super} are given in the following proposition.
\begin{cor}\label{vitexpl}
Assume that (\ref{condfeps}), (\ref{feps}), and (\ref{momg}) hold,
that $g$ belongs to $\mathcal{S}_{s,r,b}(C_1)$ defined by
\eref{super},  and that $k_n\geq n$. Assume either that
\begin{enumerate}
\item $\sum_{k\geq 1}\beta_1(k)<+\infty$
\item
or $\delta=0$, $\gamma>1$ in \eref{condfeps} and $\sum_{k\geq
1}\tau_1(k)<+\infty$
\item
or $\delta>0$ in \eref{condfeps} and $\sum_{k\geq 1}\tau_1(k)<+\infty$. \end{enumerate}
Then $\hat g_m^{(n)}$ defined
by \eref{tronque} satisfies
\begin{equation}\label{borneexplicit}
{\mathbb E}\|g-\hat g_m^{(n)}\|^2 \leq \frac{C_1}{2\pi} (m^2\pi^2+1)^{-s}\exp\{-2b
\pi^{r}m^{r}\} + \frac{2\lambda_1(f_\varepsilon,\kappa_0)\Gamma(m)
}{n}+\frac{C_2}{n}\Gamma(m)o_m(1),
\end{equation}
where $C_1$ and $C_2$ are finite constants. The constant $C_2$ depends on $\sum_{k\geq
1}\beta_1(k)$ (respectively on $\sum_{k\geq 1}\tau_1(k)$).
\end{cor}

If $\gamma=1$ when $\delta=0$, then the bound \ref{borneexplicit} becomes
\begin{equation}
{\mathbb E}\|g-\hat g_m^{(n)}\|^2 \leq \frac{C_1}{2\pi} (m^2\pi^2+1)^{-s}\exp\{-2b
\pi^{r}m^{r}\} + \frac{(2+C_2)\lambda_1(f_\varepsilon,\kappa_0)\Gamma(m)
}{n},
\end{equation}
with  $C_2$ depending on $\sum_{k\geq
1}\beta_1(k)$ (respectively on $\sum_{k\geq 1}\tau_1(k)$).

The rate of convergence of $\hat g^{(n)}_{\breve m}$  is the same as
the rate for density deconvolution for i.i.d. sequences. Our context here encompasses the particular case considered
by van Es {\it et al.}~\citeyear{vanes05}.

Table~1 below gives a summary of these rates obtained when
minimizing the right hand of (\ref{borneexplicit}). The $\breve m_g$
denotes the corresponding minimizer (see \ref{oraclepar}).

{\small
\begin{table}[ptbh]
\caption{Choice of ${\breve m}_g$ and corresponding rates under
Assumptions \eref{condfeps} and \eref{super}.}\label{rates}
\begin{center}
{
\begin{tabular}{clcc}\cline{3-4}\cline{3-4}
\multicolumn{2}{c}{} &\multicolumn{2}{c}{$f_\varepsilon$} \\
\multicolumn{2}{c}{} & $\delta=0$ & $ \delta>0$ \\
\multicolumn{2}{c}{} & ordinary smooth & supersmooth \\\hline
\multirow{8}{.2cm}{\\\vfill\null $g$} & $\;$ & $\;$ & $\;$ \\
& $\begin{array}{l}
  r=0\\
  \small{\mbox{Sobolev}(s)}
\end{array}$ &
$\begin{array}{l}
  \pi {\breve m}_g=O(n^{1/(2s+2\gamma +1)})\\
  \mbox{rate}=O(n^{-2s/(2s+2\gamma+1)})
\end{array}$  &
$\begin{array}{l}
  \pi {\breve m}_g=[\ln(n)/(2\mu+1)]^{1/\delta}\\
  \mbox{rate}=O( (\ln(n))^{-2s/\delta})
\end{array}$ \\
\cline{2-4}
& $\begin{array}{l}
  r>0\\
  \mathcal{C}^\infty
\end{array}$ &
$\begin{array}{l} \\
  \pi {\breve m}_g=\left[{\ln(n)/2b}\right]^{1/r} \\
  \mbox{ rate}= \displaystyle  O\left(\frac{\ln(n)^{(2\gamma+1)/r}}n\right)
  \;\; \end{array}$ &
$\begin{array}{c}
  {\breve m}_g  \mbox{ solution of } \\
  {{\breve m}_g}^{2s+2\gamma+1-r}\exp\{2\mu
 (\pi {\breve{m}_g})^\delta+2b \pi^r {{\breve m}_g}^r\}\\
  \qquad= O(n)
\end{array}$\\
\hline
\end{tabular}}
\end{center}
\end{table}
}
When $r>0, \delta>0$ the value of  ${\breve m}_g$ is not explicitly
given. It is obtained as the solution of the equation
\begin{equation*} {{\breve m}_g}^{2s+2\gamma+1-r}\exp\{2\mu
 (\pi {\breve{m}_g})^\delta+2b \pi^r {{\breve m}_g}^r\}= O(n).
 \end{equation*}
Consequently, the rate of $\hat g^{(n)}_{\breve m_g}$ is not easy to give explicitly
and depends on the ratio $r/\delta$. If $r/\delta$ or $\delta/r$
belongs to  $ ]k/(k+1);(k+1)/(k+2)]$ with $k$ integer, the rate of
convergence can be expressed as a function of $k$. We refer to Comte
\textit{et al.}~\citeyear{CRT1} for further discussions about those
rates. We refer to Lacour~\citeyear{LAC} for explicit formulae for the rates in the special case
$r>0$ and $\delta>0$.

\subsection{Adaptive bound}

Theorem \ref{genepenmelarch} below gives a general bound which holds
under weak dependency conditions, for $\varepsilon$ being either
ordinary or super smooth.

For $a>1$, let pen$(m)$ be defined by
\begin{equation}\label{penalitess}
\displaystyle {\rm pen}(m)=\left\lbrace\begin{array}{l}
\displaystyle 192a\frac{\Delta(m)}{n}\,\mbox{
if } 0\leq \delta< 1/3,\\
\displaystyle 64a\lambda_3\frac{ \Delta(m)\,
m^{\min((3\delta/2-1/2)_+,\delta))}}{n}\,\mbox{ if }\delta \geq
1/3,
\end{array}\right.
\end{equation}
where $\Delta(m)$ is defined by (\ref{Delta1}). 
The constant $\lambda_1(f_\varepsilon,\kappa_0)$ is defined in
\eref{lambda1} and
\begin{equation}\label{lambda2} \lambda_3=
1+\frac{32\mu\pi^\delta}{\lambda_1(f_\varepsilon,\kappa_0')}\left((\sqrt{2}+8) \|f_\varepsilon\|_{\infty}\kappa_0^{-1}\sqrt{\lambda_1(f_\varepsilon,\kappa_0)}\ind_{0\leq
\delta\leq 1}+2\lambda_1(f_\varepsilon,\kappa_0)\ind_{\delta>1}
\right).
\end{equation}
The important point here is that $\lambda_3$ is known. Hence the
penalty is explicit up to a numerical multiplicative constant.
This procedure has already been practically studied for
independent sequences $(X_t)_{t\geq 1}$ and
$(\varepsilon_t)_{t\geq 1}$ in Comte {\it et
al.}~\citeyear{CRT2,CRT1}. In particular, the  practical
implementation of the penalty functions, and the calibration of
the constants have been studied in the two previously mentioned
papers. Moreover, it is shown therein that the estimation
procedure is robust to various types of dependence, whether the
errors $\varepsilon_i$'s are ordinary or super smooth
(see Tables 4 and 5 in Comte {\it et al.}~\citeyear{CRT2}).\\

In order to bound up $\mathrm{pen}(m)$, we impose that
\begin{eqnarray}
\label{mn} \pi m_n\leq \left\{
\begin{array}{ll}
n^{1/(2\gamma+1)} &\mbox{ if }\delta=0\\
\displaystyle\left[\frac{\ln(n)}{2\mu}+\frac{2\gamma+1-\delta}{2\delta\mu
  }\ln\left(\frac{\ln(n)}{2\mu}\right)\right]^{1/ \delta} &\mbox{ if
} \delta>0.
\end{array}
\right.
\end{eqnarray}
Subsequently we set
\begin{eqnarray}
\label{kappaa}
C_a=\max(\kappa_a^2,2\kappa_a) \mbox{ where } \kappa_a=(a+1)/(a-1).
\end{eqnarray}

\begin{thm}\label{genepenmelarch}
Assume that $f_\varepsilon$ satisfies \eref{condfeps}  and \ref{feps},  that $g$ satisfies \eref{momg}, and that $m_n$
satisfies \eref{mn}. Let $\mathrm{pen}(m)$ be defined by
\eref{penalitess}. Consider the collection of estimators $\hat
g_m^{(n)}$ defined by \eref{tronque} with $k_n\geq n$ and $1\leq
m\leq m_n$. Let $\beta_{\infty}$ and $\tau_{\infty}$ be defined as in (\ref{betainfini}) and (\ref{tauinfini}) respectively.
Assume either that
\begin{enumerate}
\item $\beta_{\infty}(k)=O(k^{-(1+ \theta)})$ for some $\theta
>3$
\item or $\delta=0$, $\gamma\geq 3/2$  in \eref{condfeps} and  $\tau_{\infty}(k)=
O(k^{-(1+ \theta)})$ for some $\theta
>3+2/(1+2\gamma)$
\item  or $\delta>0$ in \eref{condfeps} and $\tau_{\infty}(k)=
O(k^{-(1+ \theta)})$ for some $\theta
>3$.
\end{enumerate}
Then the  estimator  $\tilde g = \hat
g_{\hat m}^{(n)}$ defined by \eref{estitronc} satisfies
\begin{equation}\label{adaptatif}
{\mathbb E}(\|g-\tilde g\|^2) \leq C_a\inf_{m\in \{1, \cdots,
m_n\}}\Big [\|g-g_m\|^2+ {\rm pen}(m)+\frac{m^2(M_2+1)}{n} \Big ] +
\frac{\overline{C}}{n},
\end{equation} where $C_a$ is defined in \eref{kappaa} and
 $\overline{C}$ is a constant depending on $f_{\varepsilon}$, $a$,
and  $\sum_{k\geq 1}\beta_{\infty}(k)$ (respectively on $\sum_{k\geq 1}\tau_{\infty}(k)$).

\end{thm}

\begin{rem}{\rm
In case (2), when $\delta=0$ in \eref{condfeps}, the condition on
$\theta$  is weaker as $\gamma$ increases and $f_\varepsilon$ gets
smoother. }
\end{rem}

The estimator $\tilde g$ is adaptive in the sense that it is purely
data-driven. This is due to the fact that pen$(.)$ is
explicitly known. In particular, its construction does not require any prior smoothness
knowledge on the unknown density $g$ and does not use  the
dependency coefficients. This point is important since all quantities involving dependency coefficients are usually not tractable in practice.

  The main result in Theorem \ref{genepenmelarch} shows that
the MISE of $\tilde g$ automatically achieves the best
squared-bias variance compromise (possibly up to some logarithmic
factor). Consequently, it achieves the best rate among the rates
of the $\hat g_m^{(n)}$, even from a non-asymptotical point of
view. This last point is of most importance since the $m$ selected
in practice are small and far away from asymptotic. For practical
illustration of this point in the case of density deconvolution of
i.i.d. variables, we refer to
 Comte {\it et
al.}~\citeyear{CRT2,CRT1}. Another important point is that, if we consider the
asymptotic trade-off, then  the rates given in
Table \ref{rates} are automatically reached in most cases by the adaptive estimator $\tilde g$.
 Only in the case
$\delta>1/3$ and $r>0$, a loss may occur in the rate of $\tilde g$. This comes from the additional power of $m$ in
the penalty for $\delta\geq 1/3$ with respect to the variance order $\Delta(m)$. Nevertheless, the resulting loss in the rate
has an order which is negligible compared to
the main order rate.

As a conclusion, the estimator $\tilde g$ has the rate of the i.i.d. case, with an explicit
penalty function not depending on the dependency coefficients.

\section{Proofs}

\setcounter{equation}{0}
\setcounter{lem}{0}
\setcounter{thm}{0}
\subsection{Proof of Proposition \ref{Vitss}}
The proof of Proposition \ref{Vitss} follows the same lines as
in the independent framework (see Comte \textit{ et
al.}~\citeyear{CRT1}). The main difference lies in the control of
the variance term. We keep the same notations as in Section
\ref{nonpen}. According to \eref{tronque}, for any given $ m$
belonging to $\{1,\cdots,m_n\}$, $\hat g_{m}^{(n)}$ satisfies,
$\gamma_n(\hat g_{m}^{(n)} )-\gamma_n(g_{m}^{(n)})\leq 0.$ For a
random variable $T$ with density $f_T$, and any function $\psi$ such
that $\psi(T)$ is integrable, set $\nu_{n,T}(\psi)=n^{-1}\sum_{i=1}^n [\psi(T_i)-\langle \psi,
f_T\rangle]$. In particular,
\begin{equation}\label{nu}
\nu_{n,Z}(u_t^*)=\frac 1n
\sum_{i=1}^n \left[u_t^*(Z_i)-\langle t,g\rangle
\right].\end{equation} Since
\begin{eqnarray}
\label{difgamma}
\gamma_n(t)-\gamma_n(s)=\|t-g\|^2-\|s-g\|^2-2\nu_{n,Z}(u_{t-s}^*),
\end{eqnarray}
we infer that \begin{equation}\label{premineq}\|g-\hat
g_m^{(n)}\|^2\leq \|g-g_m^{(n)}\|^2 + 2\nu_{n,Z}\left(u_{\hat
g_m^{(n)} - g_m^{(n)}}^*\right) \, . \end{equation}
 Writing that $\hat
a_{m,j}-a_{m,j}= \nu_{n,Z}(u_{\varphi_{m,j}}^*)$, we obtain that
$$
\nu_{n,Z}\left(u_{\hat g_m^{(n)}-g_m^{(n)}}^*\right)=\sum_{\vert
j\vert \leq k_n} (\hat
a_{m,j}-a_{m,j})\nu_{n,Z}(u_{\varphi_{m,j}}^*) = \sum_{\vert
j\vert\leq k_n} [\nu_{n,Z}(u_{\varphi_{m,j}}^*)]^2.$$
Consequently,
$\esp\|g-\hat g_m^{(n)}\|^2\leq  \|g-g_m^{(n)}\|^2 +
2\sum_{j\in \mathbb{Z}}\mathbb{E}[(\nu_{n,Z}(u_{\varphi_{m,j}}^*))^2].$
According to Comte \textit{et al.} \citeyear{CRT1},
\begin{eqnarray} \label{biasineq}
\|g-g_m^{(n)}\|^2= \parallel
g-g_m\parallel^2+\|g_m-g_m^{(n)}\|^2\leq \parallel g-g_m\parallel^2+
\frac{(\pi m)^2(M_2+1)}{k_n}.
\end{eqnarray}
The variance term is studied by using first that for $f\in
\mathbb{L}_1(\mathbb{R})$,
\begin{eqnarray}
\label{lemnu} \nu_{n,Z}(f^*)=\int \nu_{n,Z}(e^{ix\cdot})f(x)dx.
\end{eqnarray}
Now, we use \eref{lemnu}
and apply Parseval's formula to obtain
\begin{eqnarray}\nonumber
\mathbb{E}\Big(\sum_{j\in
\mathbb{Z}}(\nu_{n,Z}(u_{\varphi_{m,j}}^*))^2\Big)&=&\frac{1}{4\pi^2}\sum_{j\in
\mathbb{Z}}\mathbb{E}\Big(\int
\frac{\varphi_{m,j}^*(-x)}{f_\varepsilon^*(x)}\nu_{n,Z}(e^{ix\cdot})dx\Big)^2\\
\label{varineq} &=& \frac{1}{2\pi}\int_{-\pi m}^{\pi
m}\frac{\mathbb{E}\vert \nu_{n,Z}(e^{ix\cdot})\vert^2}{\vert
f^*_\varepsilon(x)\vert^2}dx.
\end{eqnarray}
Since $\nu_{n,Z}$ involves centered and stationary variables, we have
\begin{equation}
\mathbb{E}\vert \nu_{n,Z}(e^{ix\cdot})\vert^2=\mbox{Var}\vert
\nu_{n,Z}(e^{ix\cdot})\vert
=\frac{1}{n}\mbox{Var}(e^{ixZ_1})+\frac{1}{n^2}\sum_{1\leq
k\not=l\leq n}\mbox{Cov}(e^{ixZ_k},e^{ixZ_l}). \label{devpt}
\end{equation}
It follows from the structure of  the model that, for $k<l$,
$\varepsilon_l$ is independent of $(X_l,Z_k)$, so that
$\mathbb{E}(e^{ix Z_k})=f_\varepsilon^*(x)g^*(x)$ and ${\mathbb
E}(e^{ix(Z_l-Z_k)})= f_{\varepsilon}^*(x){\mathbb
E}(e^{ix(X_l-Z_k)})$. Thus, for $k<l$,
\begin{eqnarray}
\label{inegcov}
\mbox{Cov}(e^{ixZ_k},e^{ixZ_l})=f_{\varepsilon}^*(x)\mbox{Cov}(e^{ixZ_k},e^{ixX_l}).
\end{eqnarray}
>From \eref{devpt} and the stationarity of $(X_i)_{i \geq 1}$, we
obtain that
\begin{equation}\label{nunfin}
\mathbb{E}\vert \nu_{n,Z}(e^{ix\cdot})\vert^2 \leq
\frac{1}{n}+\frac{2}{n}\sum_{k=2}^n\left\vert\mbox{Cov}(e^{ixZ_1},e^{ixX_k})\right\vert\vert
f_{\varepsilon}^*(x)\vert.
\end{equation}
The first part of Proposition \ref{Vitss} follows from  the
stationarity of the $X_i$'s, and from \eref{premineq},
\eref{biasineq}, \eref{varineq} and \eref{nunfin}.

The proof of $R_m \leq \min (R_{m, \beta}, R_{m, \tau})$, where
$R_{m, \beta}$ and $R_{m, \tau}$ are defined in Proposition
\ref{Vitss}, comes from the inequalities (\ref{covineg}) in
Section \ref{couplage}. Hence we get the result.$\Box$

\subsection{Proof of Corollary \ref{vitexpl}}

According to Butucea and Tsybakov
\citeyear{ButTsyb}, under \eref{condfeps}, we have
$$
\lambda_1(f_\varepsilon,\kappa_0')\Gamma(m)(1+o_m(1))\leq
\Delta(m)\leq \lambda_1(f_\varepsilon,\kappa_0)\Gamma(m)(1+o_m(1))
\quad \text{as $m\rightarrow\infty$, where }$$
\begin{equation}\label{gammadem}
\Gamma(m)= (1+(\pi m)^2)^{\gamma}(\pi m)^{1-\delta}
\exp\big\{2\mu(\pi m)^{\delta}\big\},
\end{equation}
where $\lambda_1$ is defined in \eref{lambda1}.
In the same way
\begin{multline*}
\overline{\lambda_1}(f_\varepsilon,\kappa_0')\overline{\Gamma}(m)
(1+o_m(1))\leq
\Delta_{1/2}(m)\leq \overline{\lambda_1}(f_\varepsilon,\kappa_0)\overline{\Gamma}(m)
(1+o_m(1))
\quad \text{as $m\rightarrow\infty$,}\end{multline*}
where
\begin{eqnarray*}
\overline{\Gamma}(m)&=&(1+(\pi m)^2)^{\gamma/2}(\pi
m)^{1-\delta}\exp(\mu (\pi m)^\delta)\\
\overline{\lambda_1}(f_\varepsilon,\kappa_0)&=&\big[\kappa_0^{2}\pi (\mbox{1}\!\mbox{I}_{\{\delta=0\}} + \mu\delta
\mbox{1}\!\mbox{I}_{\{\delta>0\}})\big]^{-1} .\end{eqnarray*}
It is easy to see that $\Delta_{1/2}(m)\leq \sqrt{m\Delta(m)}$ and hence
$\Delta_{1/2}(m)=\Gamma(m)o_m(1)$. Now,  as soon as $\gamma>1$
when $\delta=0$, $m\Delta_{1/2}(m)=\Gamma(m)o_m(1).$
Set
$m_1$ such that for $m\geq m_1$ we have
\begin{equation}
\label{encadDelta}
0.5\lambda_1(f_\varepsilon,\kappa_0')\Gamma(m)\leq\Delta(m) \leq
2\lambda_1(f_\varepsilon,\kappa_0)\Gamma(m),\end{equation}
and
\begin{equation}
\label{encadDeltab}
0.5\overline{\lambda_1}(f_\varepsilon,\kappa_0')\overline{\Gamma}(m)\leq\Delta_{1/2}(m) \leq
2\overline{\lambda_1}(f_\varepsilon,\kappa_0)\overline{\Gamma}(m).\end{equation}
If $\sum_{k\geq 1}\beta_1(k)<+\infty$, \eref{condfeps} and \eref{momg}
hold, and if $k_n\geq n$, then we have the upper bounds: for
$m\geq m_1$, $\lambda_1= \lambda_1(f_\varepsilon,\kappa_0)$ and $\overline{\lambda_1}=\overline{\lambda_1}(f_\varepsilon,\kappa_0)$,
\begin{eqnarray*}
{\mathbb E}\|g-\hat g_m^{(n)}\|^2&\leq& \nonumber\|g-g_m\|^2+
\frac{m^2(M_2+1)}{n}+
\frac{2\lambda_1\Gamma(m)}{n}
+8\overline{\lambda_1}\sum_{k\geq 1}\beta_1(k)\frac{\overline{\Gamma}(m)}{n}
\\
&\leq& \|g-g_m\|^2+
\frac{m^2(M_2+1)}{n}+
\frac{2\lambda_1\Gamma(m)}{n}
+\frac{C(\sum_{k\geq 1}\beta_1(k))\Gamma(m)}{n}o_m(1)
.
\end{eqnarray*}
In the same way, if $\sum_{k\geq 1}\tau_1(k)<+\infty$, if $\gamma>1$ when $\delta=0$,
if \eref{condfeps} and \eref{momg} hold, and if $k_n\geq n$, then we have
the upper bound: for
$m\geq m_1$,
\begin{eqnarray*}
{\mathbb E}\|g-\hat g_m^{(n)}\|^2&\leq&\nonumber \|g-g_m\|^2+
\frac{m^2(M_2+1)}{n}+
\frac{2\lambda_1\Gamma(m)}{n}+
2\pi\overline{\lambda_1}\sum_{k\geq 1}\tau_1(k)\frac{m\overline{\Gamma}(m)}{n}
\\&\leq & \|g-g_m\|^2+
\frac{m^2(M_2+1)}{n}+
\frac{2\lambda_1\Gamma(m)}{n}
+\frac{C(\sum_{k\geq 1}\tau_1(k))\Gamma(m)}{n}o_m(1)
.
\end{eqnarray*}
Since $\gamma > 1$ when $\delta=0$, the residual term
$n^{-1}m^2(M_2+1)$ is negligible with respect to the variance term.

Finally, $g_m$ being the orthogonal projection of $g$ on $S_m$,
we get  $g_m^*=g^*\mbox{1}\!\mbox{I}_{[-m\pi, m\pi]}$ and
therefore
$$
\|g-g_m\|^2 = \frac{1}{2 \pi} \|g^*-g_m^*\|^2 = \frac{1}{2 \pi}
\int_{|x|\geq \pi m} |g^*|^2(x)dx.
$$
If $g$ belongs to  the class $\mathcal{S}_{s,r,b}(C_1)$ defined in
\eref{super}, then
$$
\|g-g_m\|^2 \leq \frac{C_1}{2\pi} (m^2\pi^2+1)^{-s}\exp\{-2b
\pi^{r}m^{r}\}.
$$
 The corollary is proved. $\Box$

\subsection{Proof of Theorem \ref{genepenmelarch}}
By definition, $\tilde g$ satisfies that for all $m\in
\{1,\cdots,m_n\}$, $$\gamma_n(\tilde g)+\mbox{pen}(\hat m)\leq
\gamma_n(g_m)+\mbox{pen}(m).$$ Therefore, by using \eref{difgamma}
we get 
\begin{eqnarray*}
\|\tilde g-g\|^2\leq \|g_m^{(n)}-g\|^2+2\nu_{n,Z}(u_{\tilde
g-g_m^{(n)}}^*)+\mbox{pen}(m)-\mbox{pen}(\hat m),
\end{eqnarray*}
where  $\nu_{n,Z}$ is defined in \eref{nu}.
If $t=t_1+t_2$ with $t_1$ in $S_m^{(n)}$ and $t_2$ in
$S_{m'}^{(n)}$, $t^*$ has its support in $[-\pi {\max(m,m')}, \pi
{\max(m,m')}]$ and  $t$ belongs to $S_{\max(m,m')}^{(n)}$. Set
 $B_{m, m'}(0,1)=\{t\in
S_{\max(m,m')}^{(n)} \;/\; \|t\|=1\}$ and write  $$|\nu_{n,Z}(u_{\tilde
g-g_m^{(n)}}^*) |\leq \|\tilde g-g_m^{(n)}\|\sup_{t\in
  B_{m,\hat m}(0,1)}|\nu_{n,Z}(u_t^*)|.$$
Using that $2uv \leq a^{-1}u^2+av^2$ for any $a>1$, leads to
\begin{eqnarray*} \|\tilde g-g\|^2 &\leq& \|g_m^{(n)}
-g\|^2 + a^{-1}\|\tilde g-g_m^{(n)}\|^2  + a\sup_{t\in B_{m,\hat
m}(0,1)}(\nu_{n,Z}(u_t^*))^2+ {\rm pen}(m)- {\rm pen}(\hat m).
\end{eqnarray*}

\noindent{\bf Proof in the $\beta$-mixing case.} \\
We use the coupling
methods recalled in Section \ref{couplage} to build approximating variables for the $W_i=(Z_i,X_i)$'s. More
precisely, we build variables $W_i^{\star}$ such that if
$n=2p_nq_n+r_n$, $0\leq r_n<q_n$, and $\ell =0, \cdots, p_n-1$
$$E_{\ell}=(W_{2 \ell q_n+1},...,W_{(2
  \ell+1)q_n}),~~~ F_{\ell} =(W_{(2 \ell+1)q_n+1},...,W_{(2 \ell +2)q_n}),$$
$$E_{\ell}^{\star}=(W_{2 \ell q_n+1}^{\star},...,W_{(2
  \ell+1)q_n}^{\star}),~~~ F_{\ell}^{\star} =(W_{(2 \ell
  +1)q_n+1}^{\star},...,W_{(2 \ell +2)q_n}^{\star}).$$
The variables $E_{\ell}^\star$ and $F_{\ell}^\star$ are such that
\begin{itemize}
\item[-] $E_{\ell}^{\star}$ and $E_{\ell}$ are identically
distributed. $F_{\ell}^{\star}$ and $F_{\ell}$ are identically
distributed.  \item[-] $\mathbb{P}(E_{\ell} \not=
E_{\ell}^*)\leq \beta_{\infty}(q_n)$ \mbox{ and } $\mathbb{P}(F_{\ell} \not=
F_{\ell}^*)\leq \beta_{\infty}(q_n)$, \item[-] $E_{\ell}^{\star}$
and ${\mathcal M}_0\vee \sigma(E_0,E_1,...,E_{ \ell
-1},E_0^{\star},E_1^{\star} ,\cdots,E_{\ell-1}^{\star})$ are
independent, and therefore independent of ${\mathcal
M}_{(\ell-1)q_n}$ and the same holds for the blocks
$F_{\ell}^{\star}$.
\end{itemize}
For the sake of simplicity we assume that $r_n=0$. We denote by
$(Z^{\star}_i, X_i^{\star})=W_i^\star$ the new couple of variables.
We start from
\begin{equation}
\parallel \tilde g-g\parallel^2\leq \kappa_a^2\parallel
g_m^{(n)}-g\parallel^2+a\kappa_a\sup_{t\in B_{m,\hat
m}(0,1)}\vert\nu_{n,Z}(u_t^*)\vert^2+\kappa_a(\mathrm{pen}(m)
-\mathrm{pen}(\hat{m}))\label{basedec},
\end{equation}
where $\kappa_a$ is defined in \eref{kappaa}. Using the notation \eref{nu},
we denote by $\nu_{n,Z}^{\star}(u_t^*)$ the empirical contrast
computed on the $Z_i^{\star}$. Then we write
\begin{eqnarray*}
 \|\tilde g-g\|^2  &\leq & \kappa_a^2 \|g-g_m^{(n)}\|^2 + 2a\kappa_a
\sup_{t\in B_{m,\hat
m}(0,1)}\vert\nu_{n,Z}^\star(u_t^*)\vert^2+\kappa_a({\rm pen}(m)- {\rm pen}(\hat m))\\
&&+2a\kappa_a \sup_{t\in B_{m,\hat
m}(0,1)}\vert\nu_{n,Z}^\star(u_t^*)-\nu_{n,Z}(u_t^*)\vert^2.
\end{eqnarray*}
Set
\begin{equation}\label{Wgstar} T_n^{\star }(m,m'):=\big[\sup_{t\in B_{m, m'}(0,1)}
 |\nu_{n,Z}^{\star}(t)|^2-p(m, m')\big]_+.
\end{equation}
Hence
\begin{eqnarray}
\|\tilde g-g\|^2 \nonumber &\leq & \kappa_a^2
 \|g-g_m^{(n)}\|^2+
 2 a\kappa_a T_n^{\star }(m,\hat m)
+\kappa_a\left(2a p(m,\hat m) +{\rm pen}(m)- {\rm pen}(\hat
m)\right)\nonumber\\&&+ 2a\kappa_a\sup_{t\in B_{m,\hat m}(0,1)
}\vert\nu_{n,Z}(u_{t}^*)-\nu_{n,Z}^\star(u_t^*)\vert^2\nonumber\\\nonumber
&\leq & \kappa_a^2
 \|g-g_m^{(n)}\|^2 + 2\kappa_a{\rm pen}(m)+
2a\kappa_a\sup_{t\in B_{m,\hat m}(0,1)
}\vert\nu_{n,Z}(u_{t}^*)-\nu_{n,Z}^\star(u_t^*)\vert^2 \\
\label{ineg1} && +
 2 a\kappa_a T_n^{\star }(m,\hat m)
\end{eqnarray} where pen$(m)$ is chosen such that
\begin{eqnarray}\label{pmmp}
2ap(m,m')\leq {\rm pen}(m)+ {\rm pen}(m').\end{eqnarray}
Now write 
\begin{eqnarray*}
\nu_{n,Z}(u_{t}^*)-\nu_{n,Z}^\star(u_t^*)&=&\frac{1}{2\pi}\frac{1}{n}\sum_{k=1}^n
\int
[e^{ixZ_k}-e^{ixZ_k^\star}]\frac{t^*(-x)}{f_\varepsilon^*(x)}dx\\ &=& \frac{1}{2\pi}\int
[\nu_{n,Z}(e^{ix\cdot})-\nu_{n,Z}^\star(e^{ix\cdot})]\frac{t^*(-x)}{f_\varepsilon^*(x)}dx.
\end{eqnarray*}
Consequently,
\begin{multline}
\label{b1}\mathbb{E}\Big[\sup_{t\in B_{m,\hat m}(0,1)
}\vert\nu_{n,Z}(u_{t}^*)-\nu_{n,Z}^\star(u_t^*)\vert^2\Big]\leq
\int_{-\pi {m_n}}^{\pi {m_n}} \mathbb{E}[\vert
\nu_{n,Z}(e^{ix\cdot})-\nu_{n,Z}^\star(e^{ix\cdot})\vert^2]\frac{1}{\vert
f_\varepsilon^*(x)\vert^2}dx.
\end{multline}
Since
\begin{eqnarray*}
\mathbb{E}[\vert
\nu_{n,Z}(e^{ix\cdot})-\nu_{n,Z}^\star(e^{ix\cdot})\vert^2]&=&\mathbb{E}[\vert
\nu_{n,Z}(e^{ix\cdot})-\nu_{n,Z}^\star(e^{ix\cdot})\ind_{Z_k\not=
Z_k^\star}\vert^2]\\ &\leq&
4\mathbb{E}\Big[\frac{1}{n}\sum_{k=1}^n\vert\ind_{Z_k\not=Z_k^\star}\vert^2\Big]
\leq 4\beta_{\infty}(q_n),
\end{eqnarray*}
we obtain that
\begin{eqnarray}
\mathbb{E}\Big[\sup_{t\in B_{m,\hat m}(0,1)
}\vert\nu_{n,Z}(u_{t}^*)-\nu_{n,Z}^\star(u_t^*)\vert^2\Big]&\leq&
4\beta_{\infty}(q_n)\Delta(m_n).\label{ineg2}\end{eqnarray} By
gathering \eref{ineg1} and \eref{ineg2} we get 
\begin{eqnarray*}
\qquad\mathbb{E}\|\tilde g-g\|^2 \leq\kappa_a^2 \|g-g_m^{(n)}\|^2+
2a\kappa_a\!\!\!\sum_{m'=1}^{m_n}
\mathbb{E}\big[T_n^{\star }(m,m')\big] + 2\kappa_a{\rm
pen}(m)+ 2a\kappa_a\beta_{\infty}(q_n) \Delta(m_n).
\end{eqnarray*}
Therefore we infer that, for all $m\in \{1,\cdots,m_n\}$,
\begin{eqnarray}
\label{b3}
\mathbb{E}\|g-\tilde g\|^2 \leq
C_a\left[
 \|g-g_m^{(n)}\|^2 +   {\rm pen}(m)\right] + 2a\kappa_a(C_1+C_2)/n,\end{eqnarray} provided that
\begin{equation}\label{betaq}\Delta(m_n)\beta_\infty(q_n)\leq C_1/n\quad
\mbox{ and }\quad
 \sum_{m'=1}^{m_n}
  \mathbb{E}(T_n^{\star }(m, m'))\leq C_2/n.\end{equation} Using \eref{encadDelta},
  we conclude that the first part of
(\ref{betaq}) is fulfilled as soon as
\begin{equation}\label{contrainte} {m_n}^{2\gamma+1-\delta}
\exp\{2\mu\pi^\delta {m_n}^\delta\} \beta_\infty(q_n)\leq C_1'/n.
\end{equation}
In order to ensure that our estimators converge, we only consider
models with
bounded penalty, and therefore (\ref{contrainte}) requires that
$\beta_{\infty}(q_n)\leq C'_1/n^2$. For  $q_n=[n^c]$ and
$\beta_{\infty}(k)=O(n^{-1-\theta})$, we obtain the condition
$n^{-c(1+\theta)}=O(n^{-2})$. If $\theta>3$, one can find  $c\in
]0,1/2[$, such that this condition is satisfied. Consequently,
\eref{contrainte} holds.

To prove the second part of (\ref{betaq}), we split
$T_n^{\star }(m,m')$ into two terms
$$T_n^{\star }(m,m')=(T_{n,1}^{\star }(m,m')+T_{n,2}^{\star }(m,m'))/2,$$
where, for $k=1,2$,
\begin{multline}
\label{Rk}
T_{n,k}^{\star }(m,m')=  \Big[\sup_{t\in B_{m, m'}(0,1)}
    \big|\frac{1}{p_nq_n}\sum_{\ell=1}^{p_n}\sum_{i=1}^{q_n}\big(u^*_t(Z^\star_{(2\ell+k-1)
      q_n+i})-\langle t, g\rangle \big)\big|^2-p_k(m,
    m')\Big]_+.
\end{multline}
We only study $T_{n,1}^{\star }(m,m')$ and conclude for
$T_{n,2}^{\star }(m,m')$ analogously. The study of
$T_{n,1}^{\star }(m,m')$ consists in applying a
concentration inequality to $\nu_{n,1}^{\star}(t)$ defined by
\begin{eqnarray}
\label{nustar} \nu_{n,1}^{\star}(t)&=&
\frac{1}{p_nq_n}\sum_{\ell=1}^{p_n}\sum_{i=1}^{q_n}\left(u^*_t(Z^\star_{2\ell
      q_n+i})-\langle t, g\rangle\right)=\frac{1}{p_n}\sum_{\ell
  =1}^{p_n}\nu_{q_n,\ell}^\star(u_t^*).
\end{eqnarray}
The random variable $\nu_{n,1}^{\star}(u_t^*)$ is considered as the sum of the
$p_n$ independent random variables $\nu_{q_n,\ell}^{\star}(t)$
defined as
\begin{eqnarray}\label{nuqn}
\nu_{q_n,\ell}^{\star}(u_t^*)=(1/q_n)\sum_{j=1}^{q_n}u_t^*(Z^{\star}_{2\ell
q_n +j})-\langle t, g\rangle.\end{eqnarray} Let  $m^*=\max(m,m')$.
Let $M_1^\star(m^*)$, $v^\star(m^*)$ and $H^\star(m^*)$ be some
terms such that $\sup_{t \in
  B_{m,m'}(0,1)}\parallel\nu_{q_n,\ell}^{\star}(u_t^*)\parallel_\infty\leq
M_1^\star(m^*)$,
$\sup_{t \in B_{m,m'}(0,1)} \mbox{Var}( \nu_{q_n,\ell}^{\star}(u_t^*) ) \leq
v^\star(m)$ and  lastly $\mathbb{E}(\sup_{t \in
  B_{m,m'}(0,1)}\vert\nu_{n,1}^{\star}(u_t^*)\vert )\leq H^\star(m^*)$.
According to Lemma \ref{l4}  we take
\begin{eqnarray*}
(H^{\star}(m^*))^2=\frac{2\Delta(m^*)}{n},
\;M_1^\star(m^*)=\sqrt{\Delta(m^*)} \mbox{ and
}v^\star(m^*)=\frac{2\sqrt{\Delta_2(m^*, f_Z)}}{2\pi q_n},
\end{eqnarray*}
where
\begin{equation}\label{Delta2}
\Delta_2(m, f_Z)=  \int_{-\pi m}^{\pi m} \int_{-\pi m}^{\pi
m}\frac{\vert f_Z^*(x-y)\vert^2}{\vert
f_\varepsilon^*(x)f_\varepsilon^*(y)\vert^2}dxdy.
 \end{equation}
 From the definition of $T_{n,1}^{\star }(m,m')$, by taking
$p_1(m,m')=2(1+2\xi^2)(H^\star)^2(m^*)$, we get 
\begin{equation}
\mathbb{E}(T_{n,1}^{\star }(m,m'))\leq\mathbb{E}\big[\sup_{t\in
B_{m,m'}(0,1)}\vert\nu_{n,1}^\star(u_t^*)
  -2(1+2\xi^2)(H^\star)^2(m^*)\big]_+.\label{bWn}\end{equation}
According to the condition \eref{pmmp}, we thus take
\begin{eqnarray}
\nonumber
\mbox{pen}(m)&=&4ap(m,m) =4a(2p_1(m,m)+2p_2(m,m))
=16ap_1(m,m)
\\&=&32a(1+2\xi^2)\big(2n^{-1}\Delta(m)\big)=64a(1+2\xi^2)n^{-1}\Delta(m).\label{pen}
\end{eqnarray}
where $\xi^2$ is suitably chosen. Set  $m_2$ and
$m_3$ as defined in Lemma \ref{l4}, and set $m_1$  such that for $m^*\geq
m_1$, $\Delta(m^*)$ satisfies \eref{encadDelta}.
Take $m_0=m_1\vee m_2\vee m_3$. We split the sum over $m'$ in two
parts and write
\begin{equation}\label{m0}
\sum_{m'=1}^{m_n}\mathbb{E}(T_{n,1}^{\star }(m,m'))=\sum_{m'|m^*\leq
m_0}\mathbb{E}(T_{n,1}^{\star }(m,m'))+ \sum_{m'|m^*\geq
m_0}\mathbb{E}(T_{n,1}^{\star }(m,m')).\end{equation} By
applying Lemma \ref{Concent}, we get
$\mathbb{E}(T_{n,1}^{\star }(m,m'))\leq K[I({m^*})+II(m^*)],$
where
$$
I(m^*)=\frac{\sqrt{\Delta_2(m^*,f_Z)}}{p_n}
\exp\left\lbrace- 2K_1\xi^2\frac{\Delta(m^*)}{v^\star(m^*)}\right\rbrace, \;
 II(m^*)= \frac{\Delta(m^*)}{p_n^2}
\exp\left\{-2K_1\xi C(\xi)\sqrt{\frac n{q_n}} \right\}.
$$
When $m^*\leq m_0$, with $m_0$ finite, we get that, for all $m\in
\{1,\cdots,m_n\}$,
\begin{eqnarray*}
\sum_{m'|m^*\leq m_0}\mathbb{E}(R_{n,1}^{\star }(m,m'))\leq
\frac{C(m_0)}{n}.
\end{eqnarray*}
We now come to the sum over $m'$ such that $m^*\geq m_0$.
It follows from Comte {\it et al.}~\citeyear{CRT1} that
\begin{eqnarray}\label{vstar} v^\star(m^*)=\frac{2\sqrt{\Delta_2(m^*,f_Z)}}{2\pi q_n}&\leq
&2\lambda_2^\star(f_\varepsilon,\kappa_0)
\frac{\Gamma_2(m^*)}{q_n},
\end{eqnarray}
with \begin{equation} \label{lambda2star}
\lambda_2^\star(f_\varepsilon,\kappa_0)= \kappa_0^{-1} \sqrt{2\pi \lambda_1}\|f_{\varepsilon^*}\| \1_{\delta\leq 1} + \1_{\delta>1}
\end{equation} where $\lambda_1=\lambda_1(f_\varepsilon,\kappa_0)$ is defined in \eref{lambda1} and
\begin{equation}\label{Gamma2} \qquad\Gamma_2(m)=(1+(\pi
m)^2)^{\gamma}(\pi
m)^{\min((1/2-\delta/2),(1-\delta))}\exp(2\mu(\pi m)^\delta)=(\pi m)^{-(1/2-\delta/2)_+}\Gamma(m).
\end{equation}
By combining the left hand-side of \eref{encadDelta} and
\eref{vstar}, we get that, for $m^*\geq m_0$,
\begin{eqnarray*}
&&I(m^*)\leq\frac{\lambda_2^\star(f_\varepsilon,\kappa_0)\Gamma_2(m^*)}{n}
\exp\left\lbrace -\frac{K_1\xi^2\lambda_1(f_\varepsilon,\kappa_0')}
{2\lambda_2^\star(f_\varepsilon,\kappa_0)}(\pi m^*)^{(1/2-\delta/2)_+}\right\rbrace\\
\mbox{ and }
&&II(m^*)\leq\frac{\Delta(m^*)q_n^2}{n^2}
\exp\left\lbrace-\frac{2K_1\xi
C(\xi)}{7}\frac{\sqrt{n}}{q_n}\right\rbrace.
\end{eqnarray*}

\noindent $\bullet$ Study of $\sum_{m'|m^*\geq m_0}II(m^*)$.
According to the choices for $v^\star(m^*)$, $(H^\star(m^*))^2$
and $M_1^\star(m^*)$, we have
\begin{eqnarray*}
\sum_{m'|m^*\geq m_0}II(m^*)&\leq&
\sum_{m'\in\{1,\cdots,m_n\}}\frac{\Delta(m^*)q_n^2}{n^2}
\exp\left\lbrace-\frac{2K_1\xi
C(\xi)}{7}\frac{\sqrt{n}}{q_n}\right\rbrace\\&=&
O\left[m_n \exp\left\lbrace -\frac{2K_1\xi
C(\xi)}{7}
 \frac{ \sqrt{n}}{q_n}\right\rbrace\frac{\Delta(m_n)q_n^2}{n^2}\right].
\end{eqnarray*}
Since  $\Delta(m_n)/n$ is bounded, then $q_n=[n^c]$ with $c$ in
$]0,1/2[$ ensures that
\begin{eqnarray}
\label{inegm} \sum_{m'=1}^{m_n} m_n \exp\left\lbrace
-\frac{2K_1\xi C(\xi)}{7}
 \frac{ \sqrt{n}}{q_n}\right\rbrace\frac{\Delta(m_n)q_n^2}{n^2}\leq \frac{C}{n}.
\end{eqnarray}
Consequently
\begin{eqnarray}
\label{inegm2}
\sum_{m'|m^*\geq m_0}II^\star(m^*)\leq \frac{C}{n}.\end{eqnarray}

\noindent $\bullet$ Study of $\sum_{m'|m^*\geq m_0}I(m^*)$. Denote
by $\psi=2\gamma+ \min(1/2-\delta/2,1-\delta)$,
$\omega=(1/2-\delta/2)_+$, and
$K'=K_1\lambda_1(f_\varepsilon,\kappa_0')/(2\lambda_2^\star(f_\varepsilon,\kappa_0)).$ For $a,b\geq
1$, we use that
\begin{eqnarray}\nonumber
&& \max(a,b)^{\psi}e^{2\mu\pi^{\delta}
\max(a,b)^{\delta}}e^{-K'\xi^2\max(a,b)^{\omega}} \leq
(a^{\psi}e^{2\mu\pi^{\delta} a^{\delta}}+b^{\psi}
e^{2\mu\pi^{\delta}
b^{\delta}})e^{-(K'\xi^2/2)(a^{\omega} + b^{\omega})}\\
&& \hspace{1.5cm}\leq    a^{\psi}e^{2\mu\pi^{\delta}
a^{\delta}}e^{-(K'\xi^2/2)a^{\omega}} e^{-(K'\xi^2/2)
b^{\omega}}+b^{\psi} e^{2\mu\pi^{\delta} b^{\delta}}
e^{-(K'\xi^2/2) b^{\omega}}.\label{eqmax}\end{eqnarray}
Consequently,
\begin{eqnarray}
&& \sum_{m'|m^*\geq m_0}I(m^*)\leq \!\!\!
\sum_{m'=1}^{m_n}\frac{\lambda_2^\star(f_\varepsilon,\kappa_0)\Gamma_2(m^*)}{n}
\exp\left\lbrace-\frac{K_1\xi^2\lambda_1(f_\varepsilon,\kappa_0')}{2\lambda_2^\star(f_\varepsilon,\kappa_0)}(\pi m^*)^{(1/2-\delta/2)_+}\right\rbrace\nonumber\\
&\leq &
\frac{2\lambda_2^\star(f_\varepsilon,\kappa_0)\Gamma_2(m)}{n}\exp\left\lbrace-\frac{K'\xi^2}{2}(\pi
m)^{(1/2-\delta/2)_+}\right\rbrace\sum_{m'=1}^{m_n} \exp\left\lbrace-\frac{K'\xi^2}{2}(\pi m')^{(1/2-\delta/2)_+}\right\rbrace\nonumber \\
&& \hspace{1cm} +
\sum_{m'=1}^{m_n}\frac{2\lambda_2^\star(f_\varepsilon,\kappa_0)\Gamma_2(m')}{n}\exp\left\lbrace
-\frac{K'\xi^2}{2}(\pi m')^{(1/2-\delta/2)_+}\right\rbrace
\label{I}.
\end{eqnarray}

\paragraph{\textbf{Case $0\leq\delta < 1/3$}} In that case, since
$\delta< (1/2-\delta/2)_+$, the choice $\xi^2=1$ ensures that
$\Gamma_2(m)\exp\{-(K'\xi^2/2)({m})^{(1/2-\delta/2)}\}$ is bounded
and thus the first term in \eref{I} is bounded by $C/n.$ Since
$1\leq m\leq m_n$ with $m_n$ such that $\Delta(m_n)/n$ is bounded, the term
$\sum_{m'=1}^{m_n}\Gamma_2(m')\exp\{-(K'/2)({m'})^{(1/2-\delta/2)}\}/n
$ is bounded by $C'/n,$ and hence $$\sum_{m'|m^*\geq
m_0}I(m^*)\leq \frac{C}{n}.$$ According to \eref{pmmp}, the result follows
by choosing $\mbox{pen}(m)=4ap(m,m)=192a\Delta(m)/n. $

\medskip
\paragraph{\textbf{Case $\delta= 1/3$}} According to the inequality
\eref{eqmax}, $\xi^2$ is such that $2\mu\pi^{\delta}(m)^{\delta} -
(K'\xi^2/2){m}^{\delta}= -2\mu (\pi {m^*})^{\delta}$ that is
\begin{eqnarray*}
\xi^2=\frac{16\mu\pi^\delta\lambda_2^\star(f_\varepsilon,\kappa_0)}{K_1\lambda_1(f_\varepsilon,\kappa_0')}.\end{eqnarray*}
Arguing as for the case $0\leq \delta<1/3$, this choice ensures
that $\sum_{m'|m^*\geq m_0}I(m^*) \leq C/n$. The result follows by
taking
$p(m,m')=2(1+2\xi^2)
\Delta(m^*)/ n,$ and
\begin{eqnarray*}
\mbox{pen}(m)&=&64a(1+2\xi^2)\frac{\Delta(m)}{n}
= 64a\left(1+\frac{32\mu\pi^\delta\lambda_2^\star(f_\varepsilon,\kappa_0)}{K_1\lambda_1(f_\varepsilon,\kappa_0')}\right)\frac{\Delta(m)}{n}.
\end{eqnarray*}
\medskip

\paragraph{\textbf {Case $\delta> 1/3$}} In that case $\delta >
(1/2-\delta/2)_+$. We choose $\xi^2$ such that
$$2\mu\pi^{\delta}({m^*})^{\delta} -
(K'\xi^2/2)({m^*})^{\omega}= -2\mu\pi^{\delta}({m^*})^{\delta}.$$ In
other words
\begin{eqnarray*}
\xi^2=\xi^2(m^*)=\frac{16\mu(\pi)^{\delta}\lambda_2^\star(f_\varepsilon,\kappa_0)}{
K_1\lambda_1(f_\varepsilon,\kappa_0')}(\pi
m^*)^{\min((3\delta/2-1/2)_+,\delta)}.
\end{eqnarray*} Hence
$\sum_{m'|m^*\geq m_0}I(m^*) \leq C/n$. The result follows by
choosing
$p(m,m')=2(1+2\xi^2(m,m'))\Delta(m)/ n,$ associated to
\begin{eqnarray*}
\mbox{pen}(m)&=&64a(1+2\xi^2(m))\frac{\Delta(m)}{n}\\
&=&64a\left(1+\frac{32\mu \pi^{\delta}\lambda_2^\star(f_\varepsilon,\kappa_0)}{
K_1\lambda_1(f_\varepsilon,\kappa_0')}(\pi
m^*)^{\min((3\delta/2-1/2)_+,\delta)}
\right)\frac{\Delta(m)}{n}\qquad\Box
\end{eqnarray*}

\medskip

\noindent{\bf Proof in the $\tau$-dependent case.} \\
We use the coupling properties recalled in Section \ref{couplage} to build approximating variables for the $W_i=(Z_i,X_i)$'s. More
precisely, we build variables $W_i^{\star}$ such that if
$n=2p_nq_n+r_n$, $0\leq r_n<q_n$, and $\ell =0, \cdots, p_n-1$
$$E_{\ell}=(W_{2 \ell q_n+1},...,W_{(2
  \ell+1)q_n}),~~~ F_{\ell} =(W_{(2 \ell+1)q_n+1},...,W_{(2 \ell +2)q_n}),$$
$$E_{\ell}^{\star}=(W_{2 \ell q_n+1}^{\star},...,W_{(2
  \ell+1)q_n}^{\star}),~~~ F_{\ell}^{\star} =(W_{(2 \ell
  +1)q_n+1}^{\star},...,W_{(2 \ell +2)q_n}^{\star}).$$
The variables $E_{\ell}^\star$ and $F_{\ell}^\star$ are such that\\
- $E_{\ell}^{\star}$ and $E_{\ell}$ are identically
distributed, $F_{\ell}^{\star}$ and $F_{\ell}$ are identically
distributed,\\ -
$\displaystyle \sum_{i=1}^{q_n}\mathbb{E}(\|W_{2\ell q_n+i}-W^\star_{2\ell q_n+i}\|_{\mathbb{R}^2})
\leq q_n\tau_{\infty}(q_n),$ $\displaystyle
 \sum_{i=1}^{q_n}\mathbb{E}(\|W_{(2\ell+1) q_n+i}-W^\star_{(2\ell+1) q_n+i}\|_{\mathbb{R}^2})
\leq q_n\tau_{\infty}(q_n),$\\
- $E_{\ell}^{\star}$
and ${\mathcal M}_0\vee \sigma(E_0,E_1,...,E_{ \ell
-1},E_0^{\star},E_1^{\star} ,\cdots,E_{\ell-1}^{\star})$ are
independent, and therefore independent of ${\mathcal
M}_{(\ell-1)q_n}$ and the same holds for the blocks
$F_{\ell}^{\star}$.\\

For the sake of simplicity we assume that $r_n=0$. We denote by
$(Z^{\star}_i, X_i^{\star})=W_i^\star$ the new couple of variables.


As for the proof in the $\beta$-mixing framework, we start from \eref{ineg1}
 with $R_n^{\star }(m,\hat m)$ defined by
\eref{Wgstar} and pen$(m)$ chosen such that \eref{pmmp} holds.
Next we use \eref{b1}
and the bound $|e^{-ixt}-e^{-ixs}|\leq \vert x\vert \vert t-s\vert$. Hence we conclude that
$$
 \sum_{i=1}^{q_n} {\mathbb E} (|e^{-i X_{2 \ell q_n+i}}-e^{-i X_{2 \ell
  q_n+i}^\star}|) \leq q_n |x|
\tau_{\mathbf{X},\infty}(q_n)
$$
It follows that
\begin{eqnarray}
\mathbb{E}\Big[\sup_{t\in B_{m,\hat m}(0,1)
}\vert\nu_{n,Z}(u_{t}^*)-\nu_{n,Z}^\star(u_t^*)\vert^2\Big]
&\leq&
\frac{1}{\pi}\int_{-\pi m_n}^{\pi m_n} \mathbb{E}\vert
\nu_{n,X}^\star(e^{ix\cdot})- \nu_{n,X}(e^{ix\cdot})\vert dx
\nonumber \\&\leq& \frac{ \tau_{\mathbf{X},\infty}(q_n)}{\pi}\int_{-\pi m_n}^{\pi m_n}
\frac{|x|}{\vert
f_\varepsilon^*(x)\vert^2}dx\nonumber \\\label{ineg2tau}
&\leq&
  \tau_{\mathbf{X},\infty}(q_n) {m_n}\Delta(m_n).
\end{eqnarray}
By gathering \eref{ineg1} and \eref{ineg2tau} we get 
\begin{eqnarray*}
\mathbb{E}\|\tilde g-g\|^2 \leq \kappa_a^2 \|g-g_m^{(n)}\|^2+
2a\kappa_a\!\!\!\sum_{m'=1}^{m_n}
\mathbb{E}\big[T_n^{\star }(m,m')\big] + 2\kappa_a{\rm
pen}(m)+ 2a\kappa_a\tau_{\infty}(q_n)m_n \Delta(m_n).
\end{eqnarray*}
Therefore we infer that, for all $m\in \{1,\cdots,m_n\}$,
\eref{b3} holds
 provided that
\begin{equation}\label{tauq}\Delta(m_n)m_n\tau_{\infty}(q_n)\leq C_1/n\quad
\mbox{ and }\quad
 \sum_{m'=1}^{m_n}
  \mathbb{E}(T_n^{\star }(m, m'))\leq
C_2/n.\end{equation} Using \eref{encadDelta},
  we conclude that the first part of
(\ref{tauq}) is fulfilled as soon as
\begin{equation}\label{contraintetau} {m_n}^{2\gamma+2-\delta}
\exp\{2\mu\pi^\delta {m_n}^\delta\} \tau_\infty(q_n)\leq C_1'/n.
\end{equation}
In order to ensure that our estimators converge, we only consider
models with
bounded penalty, that is $\Delta(m_n)=O(n)$. Therefore (\ref{contraintetau}) requires that
$m_n\tau_{\infty}(q_n)\leq C'_1/n^2$. For  $q_n=[n^c]$ and
$\tau_{\infty}(k)=O(n^{-1-\theta})$, we obtain the condition
\begin{eqnarray}
\label{condition}
m_n n^{-c(1+\theta)}=O(n^{-2}).\end{eqnarray} If $f_\varepsilon$ satisfies \eref{condfeps} with $\delta>0$, and if
$\theta>3$, one can find  $c\in
]0,1/2[$, such that \eref{condition} is satisfied.  Now, if $\delta=0$ and $\gamma\geq 3/2$ in \eref{condfeps} and if $\theta>3+2/(1+2\gamma)$,
then one can find  $c\in
]0,1/2[$, such that \eref{condition} is satisfied. These conditions ensure that
\eref{contrainte} holds.

In order to prove the second part of \eref{tauq}, we proceed as for the proof of the second part of \eref{betaq} and
split
$T_n^{\star }(m,m')$ into two terms
$$T_n^{\star }(m,m')=(T_{n,1}^{\star }(m,m')+T_{n,2}^{\star }(m,m'))/2,$$
where  the $T_{n,k}^{\star }(m,m')$'s are defined in \eref{Rk}.
We only study $T_{n,1}^{\star }(m,m')$ and conclude for
$T_{n,2}^{\star }(m,m')$ analogously. As in the $\beta$-mixing framework, the study of
$T_{n,1}^{\star }(m,m')$ consists in applying a
concentration inequality to $\nu_{n,1}^{\star}(t)$ defined in \eref{nustar}
and considered as the sum of the
$p_n$ independent random variables $\nu_{q_n,\ell}^{\star}(t)$
defined as in \eref{nuqn}. Once again, set  $m^*=\max(m,m')$, and
denote by $M_1^\star(m^*)$, $v^\star(m^*)$ and $H^\star(m^*)$ the
terms such that $\sup_{t \in
  B_{m,m'}(0,1)}\parallel\nu_{q_n,\ell}^{\star}(u_t^*)\parallel_\infty\leq
M_1^\star(m^*)$,
$\sup_{t \in B_{m,m'}(0,1)} \mbox{Var}( \nu_{q_n,\ell}^{\star}(u_t^*) ) \leq
v^\star(m)$ and  lastly $\mathbb{E}(\sup_{t \in
  B_{m,m'}(0,1)}\vert\nu_{n,1}^{\star}(u_t^*)\vert )\leq H^\star(m^*)$.
According to Lemma \ref{l4tau}, we take
\begin{eqnarray*}
(H^{\star}(m^*))^2=\frac{2\Delta(m^*)}{n},
\;M_1^\star(m^*)=\sqrt{\Delta(m^*)} \mbox{ and
}v^\star(m^*)=\frac{C_{v^*}\sqrt{\Delta_2(m^*, f_Z)}}{2\pi q_n},
\end{eqnarray*}
where $\Delta_2(m, f_Z)$ is defined in \eref{Delta2} and where
\begin{eqnarray}
\label{cvstar}
C_{v^*}=2\Big[\ind_{\delta>0}+\frac{\sqrt{2}\pi^{3/2}(2\pi)^{3/2}}{\sqrt{3}}\sum_{k\geq 1}\tau_1(k)\ind_{\delta=0}\Big].
\end{eqnarray}
 From the definition of $T_{n,1}^{\star }(m,m')$, by taking
$p_1(m,m')=2(1+2\xi^2)(H^\star)^2(m^*)$, we get 
\begin{equation}
\mathbb{E}(T_{n,1}^{\star }(m,m'))\leq\mathbb{E}\Big[\sup_{t\in
B_{m,m'}(0,1)}\vert\nu_{n,1}^\star(u_t^*)
  -2(1+2\xi^2)(H^\star)^2(m^*)\Big]_+.\label{bWntau}\end{equation}
As in the $\beta$-mixing framework we take
pen$(m)=64a\Delta(m)(1+2\xi^2)/n$
where $\xi^2$ is suitably chosen (see \eref{bWntau}). Set  $m_2$ and
$m_3$ as defined in Lemma \ref{l4tau}, and set $m_1$  such that for $m^*\geq
m_1$ \eref{encadDelta} holds.
Take $m_0=m_1\vee m_2\vee m_3$ and $K'=K_1\lambda_1(f_\varepsilon,\kappa_0')/(C_{v^*}\lambda_2^\star(f_\varepsilon,\kappa_0)).$
The end of the proof is the same as in $\beta$-mixing framework, up to possible multiplicative constants.$\Box$

\subsection{Technical lemmas}

\begin{lem}
\label{details}
\begin{eqnarray}\label{b1details}
\parallel\sum_{j\in \mathbb{Z}} |u^*_{\varphi_{m,j}}|^2\parallel_\infty \leq
\Delta(m).
\end{eqnarray}
\end{lem}
The proof of Lemma \ref{details} can be found in Comte \textit{et al.} \citeyear{CRT1}.

\begin{lem}\label{l4} Assume that $\sum_{k\geq 1}\beta_1(k)<+\infty$. Then we have
\begin{eqnarray}
\label{binfbeta}
\sup_{ t \in B_{m,m'}(0,1) }\parallel
\nu^\star_{q_n,\ell}(u_t^*)\parallel_\infty \leq
\sqrt{\Delta(m^*)}
\end{eqnarray}
Moreover, there exist $m_2$ and $m_3$ such that
\begin{eqnarray*}
&&  \mathbb{E}[\sup_{t\in
B_{m,m'}(0,1)}|\nu_{n,1}^{\star}(u_t^*)|]\leq
\sqrt{2\Delta(m^*)/n} \mbox{ for } m^*\geq m_2,\\
{\rm and } &&\sup_{t\in B_{m,m'}(0,1)} {\rm
Var}(\nu^\star_{q_n,\ell}(u_t^*)) \leq
2\sqrt{\Delta_2(m^*,f_Z)}/(2\pi q_n) \mbox{ for } m^*\geq m_3,
\end{eqnarray*}
where $\Delta(m)$ and $\Delta_2(m,f_Z)$ are defined by \eref{Delta1} and
(\ref{Delta2}).
\end{lem}

\paragraph{\textbf{Proof of Lemma \ref{l4}}}
Arguing as in Lemma \ref{details} and by using Cauchy-Schwartz
Inequality and Parseval formula, we obtain that
 the first term $\sup_{ t \in
B_{m,m'}(0,1) }\parallel
\nu^\star_{q_n,\ell}(u_t^*)\parallel_\infty$ is  bounded by
$$\sup_{ t \in B_{m,m'}(0,1) }\parallel
\nu^\star_{q_n,\ell}(u_t^*)\parallel_\infty\leq \sqrt{\sum_{j\in
\mathbb{Z}}\int \left\vert
    \frac{\varphi_{m^*,j}^*(x)}{f_\varepsilon^*(x)}
    \right\vert^2dx
  }=\sqrt{\Delta(m^*)}.$$
Next \begin{eqnarray*} {\mathbb E} \Big[
    \sup_{t\in B_{m,m'}(0,1)}
    \Big\vert\nu_{n,1}^{\star}(u_t^*)\Big\vert\Big]&=& {\mathbb E} \Big[
    \sup_{t\in B_{m,m'}(0,1)}
    \Big\vert\frac{1}{p_nq_n}\sum_{\ell=1}^{p_n}\sum_{i=1}^{q_n}u^*_t(Z^\star_{2\ell
        q_n+i})-\langle t, g\rangle
    \Big\vert\Big]\\
  &\leq & \sqrt{ \sum_{j\in
      \mathbb{Z}}\mbox{Var}(\nu_{n,1}^\star(u_{\varphi_{m^*,j}}^*))}.
      \end{eqnarray*}
  By using (\ref{varineq}) we obtain
      \begin{eqnarray*}
      \sqrt{ \sum_{j\in
      \mathbb{Z}}\mbox{Var}(\nu_{n,1}^\star(u_{\varphi_{m^*,j}}^*))}
     &=&\sqrt{\sum_{j\in {\mathbb Z}} \frac
      1{p_n^2}\sum_{\ell=1}^{p_n}
      \mbox{Var}\left(\nu^{\star}_{q_n,\ell}(u_{\varphi_{m^*,j}}^*)\right)}=\sqrt{\sum_{j\in {\mathbb Z}} \frac
      1{p_n^2}\sum_{\ell=1}^{p_n}
      \mbox{Var}\left(\nu_{q_n,\ell}(u_{\varphi_{m^*,j}}^*)\right)}\\
&=&\sqrt{\sum_{j\in {\mathbb Z}} \frac
      1{p_n}\mbox{Var}\left(\nu_{q_n,1}(u_{\varphi_{m^*,j}}^*)\right)}= \sqrt{\frac 1{2\pi p_n} \int_{-\pi m^*}^{\pi m^*} \frac{{\mathbb
E}|\nu_{q_n,1}(e^{ix.})|^2}{|f_{\varepsilon}^*(x)|^2}dx }.
\end{eqnarray*}
Now, according to (\ref{nunfin}) and (\ref{ibrb})
$${\mathbb E}|\nu_{q_n,1}(e^{ix.})|^2 \leq \frac 1{q_n} + \frac 2{q_n}\sum_{k=1}^{n-1}\beta_1(k)  |f_{\varepsilon}^*(x)|.$$
This implies that
\begin{eqnarray*}
{\mathbb E}^2\Big[\sup_{t\in B_{m,m'}(0,1)}
    \Big\vert\nu_{n,1}^{\star}(u_t^*)\Big\vert\Big]
    &\leq& \frac 1{p_n}\Big(\frac 1{q_n} \Delta(m^*)+ \frac 2{q_n}\sum_{k=1}^{n-1}\beta_1(k) \Delta_{1/2}(m^*)\Big).
\end{eqnarray*}
Since $2\sum_{k\geq 1}\beta_1(k) \Delta_{1/2}(m)\leq \Delta(m)$ for $m$ large
enough, we get that, for $m^*$ large enough,
\begin{eqnarray*}
{\mathbb E}^2\Big[\sup_{t\in B_{m,m'}(0,1)}
    \Big\vert\nu_{n,1}^{\star}(u_t^*)\Big\vert\Big]
\leq  2\Delta(m^*)/n .
\end{eqnarray*}
Now, for $t\in B_{m,m'}(0,1)$ we write
\begin{eqnarray*}
{\rm Var} \Big(\frac 1{q_n}\sum_{i=1}^{q_n} u_t^*(Z_{2\ell
q_n+i}^{\star})\Big)
&=&{\rm Var} \Big(\frac 1{q_n}\sum_{i=1}^{q_n} u_t^*(Z_i)\Big) \\
&=& \frac 1{q_n^2} \Big[\sum_{k=1}^{q_n} {\rm Var}(u_t^*(Z_k)) +
2\sum_{1\leq k <l\leq q_n} {\rm Cov}(u_t^*(Z_k), u_t^*(Z_l))\Big].
\end{eqnarray*}
According to \eref{lemnu}, \eref{inegcov} and (\ref{ibrb}) we have
\begin{eqnarray*}
|{\rm Cov}(u_t^*(Z_k), u_t^*(Z_l))|&=&\Big|\int_{-\pi m^*}^{\pi m^*}\int_{-\pi
m^*}^{\pi m^*} \frac{{\rm Cov}(e^{ixZ_k},
e^{iyZ_l})t^*(x)t^*(y)}{f_\varepsilon^*(x)f_\varepsilon^*(-y)}dxdy\Big|\\
&=& \Big|\int_{-\pi m^*}^{\pi m^*}\int_{-\pi
m^*}^{\pi m^*} \frac{f_\varepsilon^*(-y){\rm Cov}(e^{ixZ_k},
e^{iyX_l})t^*(x)t^*(y)}{f_\varepsilon^*(x)f_\varepsilon^*(-y)}dxdy\Big|\\
&\leq&\int_{-\pi m^*}^{\pi m^*}\int_{-\pi
m^*}^{\pi m^*} \frac{2\beta_1(k)|t^*(x)t^*(y)|}{|f_\varepsilon^*(x)|}dxdy.
\end{eqnarray*}
Hence,
\begin{eqnarray*}
{\rm Var} \Big(\frac 1{q_n}\sum_{i=1}^{q_n} u_t^*(Z_{2\ell
q_n+i}^{\star})\Big) &\leq & \frac 1{q_n} \Big(\int_{-\pi m^*}^{\pi m^*}\int_{-\pi m^*}^{\pi m^*} \frac{f_Z^*(u-v)
t^*(u)t^*(-v)}{f_{\varepsilon}(u)f_{\varepsilon}(-v)} dudv
\\&&\qquad+ 2\sum_{k=1}^{q_n} \beta_1(k) \int_{-\pi
m^*}^{\pi m^*}\int_{-\pi m^*}^{\pi m^*}
\Big|\frac{t^*(u)t^*(v)}{f_{\varepsilon}^*(v)}\Big|dudv\Big).
\end{eqnarray*}
Following Comte {\it et al.}~\citeyear{CRT1} and applying Parseval's formula,  the first integral is less that $\sqrt{\Delta_2(m^*, f_Z)}/2\pi$. For the second one, write
$$\int_{-\pi m^*}^{\pi m^*} \int_{-\pi m^*}^{\pi
m^*}\Big|\frac{t^*(u)t^*(v)}{f_{\varepsilon}^*(v)}\Big|dudv\leq \sqrt{2\pi
m^*} \|t^*\| \sqrt{\int |t^*(v)|^2 dv \int_{-\pi m^*}^{\pi m^*}
\frac{dv}{|f_{\varepsilon}^*(v)|^2}},$$ that is
$$
\int_{-\pi m^*}^{\pi m^*} \int_{-\pi m^*}^{\pi
m^*}\Big|\frac{t^*(u)t^*(v)}{f_{\varepsilon}^*(v)}\Big|dudv\leq
 (2\pi)^2 \sqrt{m^*\Delta(m^*)}.
$$
Using that $\gamma>1/2$ if $\delta=0$, we get that
$\sqrt{m^*\Delta(m^*)}=o_m(\sqrt{\Delta_2(m^*,f_Z)})$ and hence
the result follows for $m$ large enough. $\Box$

\begin{lem}\label{l4tau} Assume that $\sum_{k\geq 1}\tau_1(k)<+\infty$.
Assume either that
\begin{enumerate}
\item $\delta=0$, $\gamma\geq 3/2$  in \eref{condfeps}
\item or $\delta>0$ in \eref{condfeps}.
\end{enumerate}
Then we have
\begin{eqnarray}
\label{binftau}
\sup_{ t \in B_{m,m'}(0,1) }\parallel
\nu^\star_{q_n,\ell}(u_t^*)\parallel_\infty \leq
\sqrt{\Delta(m^*)}
\end{eqnarray}
Moreover, there exist $m_2$ and $m_3$ such that
\begin{eqnarray*}
&&  \mathbb{E}[\sup_{t\in
B_{m,m'}(0,1)}|\nu_{n,1}^{\star}(u_t^*)|]\leq
\sqrt{2\Delta(m^*)/n} \mbox{ for } m^*\geq m_2,\\
{\rm and } &&\sup_{t\in B_{m,m'}(0,1)} {\rm
Var}(\nu^\star_{q_n,\ell}(u_t^*)) \leq C_{v^*}\sqrt{\Delta_2(m^*,f_Z)}/(2\pi q_n) \mbox{ for } m^*\geq m_3,
\end{eqnarray*}
where $\Delta(m)$ and $\Delta_2(m,f_Z)$ are defined by \eref{Delta1} and
(\ref{Delta2}) and where $C_{v^*}$ is defined in \eref{cvstar}.
\end{lem}

\paragraph{\textbf{Proof of Lemma \ref{l4tau}}}
The proof of \eref{binftau} is the same as the proof of \eref{binfbeta}.
Next, again as for the proof of Lemma \ref{l4}
 \begin{eqnarray*} {\mathbb E} \Big[
    \sup_{t\in B_{m,m'}(0,1)}
    \Big\vert\nu_{n,1}^{\star}(u_t^*)\Big\vert\Big]
&\leq & \sqrt{ \sum_{j\in \mathbb{Z}}\mbox{Var}(\nu_{n,1}^\star(u_{\varphi_{m^*,j}}^*))}
\end{eqnarray*}
with
\begin{eqnarray*}
\sqrt{ \sum_{j\in \mathbb{Z}}\mbox{Var}(\nu_{n,1}^\star(u_{\varphi_{m^*,j}}^*))}
 &=& \sqrt{\frac 1{2\pi p_n} \int_{-\pi m^*}^{\pi m^*} \frac{{\mathbb
E}|\nu_{q_n,1}(e^{ix.})|^2}{|f_{\varepsilon}^*(x)|^2}dx }.
\end{eqnarray*}
Now, according to (\ref{nunfin}) and (\ref{evident})
$${\mathbb E}|\nu_{q_n,1}(e^{ix.})|^2 \leq \frac 1{q_n} + \frac 1{q_n}\sum_{k=1}^{n-1}\tau_1(k)|x|  |f_{\varepsilon}^*(x)|.$$
This implies that
\begin{eqnarray*}
{\mathbb E}^2\Big[\sup_{t\in B_{m,m'}(0,1)}
    \Big\vert\nu_{n,1}^{\star}(u_t^*)\Big\vert\Big]
    &\leq& \frac 1{p_n}\Big(\frac 1{q_n} \Delta(m^*)+ \frac {2\pi}{q_n}\sum_{k=1}^{n-1}\tau_1(k)m \Delta_{1/2}(m^*)\Big).
\end{eqnarray*}
Since $2\pi\sum_{k\geq 1}\tau_1(k) m\Delta_{1/2}(m)\leq \Delta(m)$ for $m$
large enough, we get that for $m^*$ large enough
\begin{eqnarray*}
{\mathbb E}^2\Big[\sup_{t\in B_{m,m'}(0,1)}
    \Big\vert\nu_{n,1}^{\star}(u_t^*)\Big\vert\Big]
\leq 2\Delta(m^*)/n.
\end{eqnarray*}
Now, for $t\in B_{m,m'}(0,1)$ we write
\begin{eqnarray*}
{\rm Var} \Big(\frac 1{q_n}\sum_{i=1}^{q_n} u_t^*(Z_{2\ell
q_n+i}^{\star})\Big)
&=&{\rm Var} \Big(\frac 1{q_n}\sum_{i=1}^{q_n} u_t^*(Z_i)\Big) \\
&=& \frac 1{q_n^2} \Big[\sum_{k=1}^{q_n} {\rm Var}(u_t^*(Z_k)) +
2\sum_{1\leq k <l\leq q_n} {\rm Cov}(u_t^*(Z_k), u_t^*(Z_l))\Big].
\end{eqnarray*}
According to \eref{lemnu}, \eref{inegcov} and (\ref{evident}) and
by applying the same arguments as for the proof of Lemma \ref{l4}
we have
\begin{eqnarray*}
|{\rm Cov}(u_t^*(Z_k), u_t^*(Z_l))|
&=& \Big|\int_{-\pi m^*}^{\pi m^*}\int_{-\pi
m^*}^{\pi m^*} \frac{f_\varepsilon^*(-y){\rm Cov}(e^{ixZ_k},
e^{iyX_l})t^*(x)t^*(y)}{f_\varepsilon^*(x)f_\varepsilon^*(-y)}dxdy\Big|\\
&\leq&\int_{-\pi m^*}^{\pi m^*}\int_{-\pi
m^*}^{\pi m^*} \frac{\vert y\vert\tau_1(k)|t^*(x)t^*(y)|}{|f_\varepsilon^*(x)|}dxdy.
\end{eqnarray*}
Hence,
\begin{eqnarray*}
{\rm Var} \Big(\frac 1{q_n}\sum_{i=1}^{q_n} u_t^*(Z_{2\ell
q_n+i}^{\star})\Big)
 &\leq & \frac 1{q_n} \Big(\int_{-\pi m^*}^{\pi m^*}\int_{-\pi m^*}^{\pi m^*} \frac{f_Z^*(u-v)
t^*(u)t^*(-v)}{f_{\varepsilon}(u)f_{\varepsilon}(-v)} dudv
\\&&\qquad+ 2\sum_{k=1}^{q_n} \tau_1(k) \int_{-\pi
m^*}^{\pi m^*}\int_{-\pi m^*}^{\pi m^*}
\Big|\frac{ u t^*(u)t^*(v)}{f_{\varepsilon}^*(v)}\Big|dudv\Big).
\end{eqnarray*}
Once again the first integral is less that $\sqrt{\Delta_2(m^*, f_Z)}/2\pi$.
For the second one, write
$$\int_{-\pi m^*}^{\pi m^*} \int_{-\pi m^*}^{\pi
m^*}\Big|\frac{u t^*(u)t^*(v)}{f_{\varepsilon}^*(v)}\Big|dudv\leq
\frac{\sqrt{2}\pi^{3/2}}{\sqrt{3}}
(m^*)^{3/2} \|t^*\| \sqrt{\int |t^*(v)|^2 dv \int_{-\pi m^*}^{\pi m^*}
\frac{dv}{|f_{\varepsilon}^*(v)|^2}},$$ that is
$$
\int_{-\pi m^*}^{\pi m^*} \int_{-\pi m^*}^{\pi
m^*}\Big|\frac{t^*(u)t^*(v)}{f_{\varepsilon}^*(v)}\Big|dudv\leq
\frac{\sqrt{2}\pi^{3/2}}{\sqrt{3}} (2\pi)^{3/2} \sqrt{(m^*)^3\Delta(m^*)}.
$$
If $\delta>0$, then $\sqrt{(m^*)^3\Delta(m^*)}=o_m\sqrt{\Delta_2(m^*,f_Z)}$.
If $\gamma >3/2$ and $\delta=0$, we get that
$\sqrt{(m^*)^3\Delta(m^*)}=o_m\sqrt{\Delta_2(m^*,f_Z)}$.
Lastly, if $\gamma=3/2$ and $\delta=0$, we get that
$\sqrt{(m^*)^3\Delta(m^*)}\leq \sqrt{\Delta_2(m^*,f_Z)}$ and
the result follows for $m$ large enough. $\Box$

\begin{lem}
\label{Concent} Let $Y_1, \dots, Y_n$ be independent random
variables and  let ${\mathcal F}$ be a countable class of
uniformly bounded measurable functions. Then for $\xi^2>0$
\begin{equation*}
 \mathbb{E}\Big[\sup_{f\in {\mathcal
F}}|\nu_{n,Y}(f)|^2-2(1+2\xi^2)H^2\Big]_+ \leq \frac
4{K_1}\left(\frac vn e^{-K_1\xi^2 \frac{nH^2}v} +
\frac{98M_1^2}{K_1n^2C^2(\xi^2)} e^{-\frac{2K_1
C(\xi)\xi}{7\sqrt{2}}\frac{nH}{M_1}}\right),
\end{equation*}
with $C(\xi)=\sqrt{1+\xi^2}-1$, $K_1=1/6$, and
$$\sup_{f\in {\mathcal F}}\|f\|_{\infty}\leq M_1, \;\;\;\;
\mathbb{E}\Big[\sup_{f\in {\mathcal F}}|\nu_{n,Y}(f)|\Big]\leq H,
\;\;\;\; \sup_{f\in {\mathcal F}}\frac{1}{n}\sum_{k=1}^n{\rm
Var}(f(Y_k)) \leq v.$$\end{lem} This inequality comes from  a concentration Inequality in
Klein and Rio \citeyear{KleinRio} and arguments that can be found in Birg\'e and Massart
\citeyear{BirgeMassart98}. Usual density arguments show that this result can be applied to the class of functions ${\mathcal F}= B_{m,m'}(0,1)$.

\medskip
\noindent {\bf Proof of Proposition \ref{taumix}}. To prove (1), let
for $t>0$,  $Y_t^*= \eta_t \sigma_t^*$. Note that the sequence
$((Y_t^*, \sigma_t^*))_{t \geq 1}$ is distributed as $((Y_t,
\sigma_t))_{t \geq 1}$ and independent of ${\mathcal M}_i= \sigma
(\sigma_j, Y_j, 0 \leq j \leq i)$. Hence, by the coupling properties
of $\tau$ (see (\ref{cber})), we have that, for $n+i \leq i_1 <
\cdots < i_l$,
$$
  \tau({\mathcal M}_i, (Y^2_{i_1}, \sigma^2_{i_1}), \ldots, (Y^2_{i_l},
  \sigma^2_{i_l})) \leq \frac{1}{l} \sum_{j=1}^l \|(Y^2_{i_j},
  \sigma^2_{i_j})- ((Y^*_{i_j})^2, (\sigma_{i_j}^*))^2\|_{{\mathbb R}^2} \leq
  \delta_n \, ,
$$
and (1) follows.

To prove (2), define the function $f_\epsilon(x)=
\ln(x)\ind_{x>\epsilon}+ 2 \ln(\epsilon)\ind_{x \leq \epsilon}$ and
the function $g_\epsilon(x)=\ln(x)-f_\epsilon(x)$. Clearly, for any
$\epsilon>0$ and any $n+i \leq i_1 < \ldots < i_l$, we have
\begin{multline}\label{biz}
\tau( {\mathcal M}_i, (Z_{i_1}, X_{i_1}), \ldots, (Z_{i_l},
  X_{i_l})) \leq 2{\mathbb
  E}(|g_\epsilon(Y_0^2)|+|g_\epsilon(\sigma_0^2)|) \\ + \tau({\mathcal
  M}_i, (f_\epsilon(Y^2_{i_1}), f_\epsilon(\sigma^2_{i_1})), \ldots, (f_\epsilon(Y^2_{i_l}),
  f_\epsilon(\sigma^2_{i_l})))
\end{multline}
For $ 0 < \epsilon <1$, the function $f_\epsilon$ is
$1/\epsilon$-Lipschitz. Hence, applying (1),
$$
\tau({\mathcal
  M}_i, (f_\epsilon(Y^2_{i_1}), f_\epsilon(\sigma^2_{i_1})), \ldots, (f_\epsilon(Y^2_{i_l}),
  f_\epsilon(\sigma^2_{i_l}))) \leq \frac{ \delta_n}{\epsilon} \, .
$$
 Since $\max(f_{\sigma^2}(x), f_{Y^2}(x))\leq C |\ln(x)|^\alpha x^{-\rho}$
in a neighborhood of 0, we infer that for small enough $\epsilon$,
$$
{\mathbb
  E}(|g_\epsilon(Y^2_0)|+|g_\epsilon(\sigma^2_0)|) \leq K_1
  \epsilon^{1-\rho} |\ln(\epsilon)|^{1 + \alpha} \, ,
$$
for $K_1$ a positive constant. From (\ref{biz}), we infer that
there exists a positive constant $K_2$ such that, for small enough
$\epsilon$,
$$
\tau( {\mathcal M}_i, (Z_{i_1}, X_{i_1}), \ldots, (Z_{i_l},
  X_{i_l})) \leq K_2 \Big ( \frac{\delta_n}{\epsilon}+\epsilon^{1-\rho} |\ln(\epsilon)|^{1 +
  \alpha}\Big) \, .
$$
The result follows by taking $\epsilon=(\delta_n)^{1/(2-
\rho)}|\ln(\delta_n)|^{-(1+\alpha)/(2-\rho)}$.

 Now, we go back to the model (\ref{archinfty}). If  $\sum_{j=1}^\infty a_j <1$, the unique stationary
solution to (\ref{archinfty}) is given by Giraitis {\it et al.}
(2000):
$$
\sigma_t^2= a +a \sum_{\ell=1}^\infty \sum_{j_1, \dots,
j_l=1}^\infty a_{j_1} \ldots a_{j_l} \eta^2_{t-j_1} \ldots
\eta^2_{t-(j_1 + \cdots + j_l)} .
$$
for any $ 1 \leq k \leq n$, let
$$
  \sigma_t^2(k,n)=a+  a\sum_{\ell=1}^{[n/k]} \sum_{j_1, \dots, j_l=1}^k
a_{j_1} \ldots a_{j_l} \eta^2_{t-j_1} \ldots \eta^2_{t-(j_1 + \cdots
+ j_l)} .
$$
Clearly $$ {\mathbb E}(|\sigma^2_n-(\sigma_n^*)^2|)\leq 2 {\mathbb
E}(|\sigma_0^2- \sigma_0^2(k,n)|)\, .$$
Now
$$
{\mathbb E}
  (|\sigma^2_{0}-\sigma^2_{0}(k,n)|) \leq \Big(\sum_{l=[n/k]+1} c^l +
  \sum_{l=1}^\infty c^{l-1} \sum_{j>k} a_j \Big)\, .
$$
This being true for any $1 \leq k \leq n$, the proof of  Proposition
\ref{taumix} is complete.

\bibliographystyle{plain}
\bibliography{biblioarch}

\end{document}